\documentclass{article}
\usepackage{amsfonts,amsmath,amssymb,amsthm,amstext,txfonts}
\usepackage{extarrows,longtable}
\usepackage{graphicx,hyperref,mathrsfs,url}
\usepackage{geometry}
\usepackage{titlesec}
\usepackage[title]{appendix}

\allowdisplaybreaks

\begin{document}
\title{The parabolic Harnack inequality on non-local Dirichlet spaces in the view of pure analysis}
\author{Guanhua Liu\footnote{Fakult\"at f\"ur Mathematik, Universit\"at Bielefeld, Postfach 100131, D-33501 Bielefeld, Germany. E-mail: gliu@math.uni-bielefeld.de\\The research is supported by Deutsche Forschungsgemeinschaft (DFG) Project 317210226, SFB 1283}}
\maketitle

\newtheorem{theorem}{Theorem}[section]
\newtheorem{proposition}[theorem]{Proposition}
\newtheorem{lemma}[theorem]{Lemma}
\newtheorem{remark}[theorem]{Remark}
\newtheorem{corollary}[theorem]{Corollary}
\numberwithin{equation}{section}

\begin{abstract}
This paper provides the general theory on parabolic Harnack inequalities (PHI, for short) for regular Dirichlet forms without killing part. We prove PHI by pure analytic methods, using both Nash and Moser approaches, and yield some important properties contained in PHI. Combining our recent result on weak Harnack inequalities, we greatly enlarge the list of equivalent characterizations of PHI.
\end{abstract}

\section{Introduction}

As is well known, Harnack inequalities are important in the partial differential equation theory, controlling the ratio of a harmonic or caloric function in a regular domain. We can refer to \cite{ka} for the early history of Harnack inequalities. Let us now focus on equivalent characterizations of the parabolic Harnack inequality (PHI). For strongly local Laplacians on metric measure spaces, \cite{bbk} and \cite{ghl} provide a complete theory on this topic, respectively in the language of probability and analysis. Since the definition of a caloric function is not the same in stochastic and analytic theories, the difference lies not only in the proof, but also the inequality itself. Till now, for non-local Dirichlet forms, there is only theory of PHI in the language of probability \cite{ckw-p}. In this paper we continue the analytic investigation of weak parabolic Harnack inequalities in \cite{L0} to obtain the analog of strong PHI's, covering more equivalent conditions of PHI under weaker assumptions.

A general proof of PHI is the key of this paper. Let us recall that there are two strategies for proving PHI. The first comes from Nash \cite{nas}, using a typical blow-up technique. It was developed by Fabes and Stroock \cite{fs} for classical local operators, and also in \cite{ckw-p} for general non-local cases (stochastically). Such a proof is complicated but contains very limited information. On the contrary, the second strategy, raised by Moser \cite{mos-p}, is clear and contains information for sub- and super-caloric functions. In this paper we complete both strategies using purely analytic methods.

Throughout this paper, we assume that $(M,d)$ is a locally compact, complete and separable metric space, and $\mu$ is a Radon measure of full support on $M$. Denote
$$\bar{R}:=2\ \mathrm{diam}(M,d)\in[0,\infty].$$
For any $x\in M$ and $r>0$, let the metric ball $B(x,r)$ be the set $\{y\in M:d(x,y)<r\}$ equipped with fixed center $x$ and radius $r$. Note that every metric ball is precompact. For any $0<\lambda<\infty$ and any ball $B=B(x,r)$, let
$$\lambda B:=B(x,\lambda r)\quad\mbox{and}\quad V(x,r):=\mu(B(x,r)).$$

The following two conditions are always assumed in the sequel:

\begin{itemize}
	\item $(\mathrm{VD})$, the \emph{volume doubling} property: there exists $C\ge 1$ such that $V(x,2R)\le CV(x,R)$ for all $x\in M$ and $R>0$. Equivalently, there exist $\alpha,C>0$ such that for all $x,y\in M$ and $0<r\le R<\infty$,
	$$\frac{V(x,R)}{V(y,r)}\le C\left(\frac{R+d(x,y)}{r}\right)^\alpha.$$
	\item $(\mathrm{RVD})$, the \emph{reverse volume doubling} property: there exists $C\ge 1$ such that $V(x,R)\ge 2V(x,C^{-1}R)$ for all $x\in M$ and $0<R<\bar{R}$. Equivalently, there exist $\alpha',C>0$ such that for all $x\in M$ and $0<r\le R<\bar{R}$,
	$$\frac{V(x,R)}{V(x,r)}\ge C^{-1}\left(\frac{R}{r}\right)^{\alpha'}.$$
\end{itemize}

Let $\mathcal{F}$ be a dense subspace of $L^2(M,\mu)$, and $\mathcal{E}$ be a symmetric non-negative semi-definite quadratic form on $\mathcal{F}$. We call $(\mathcal{E},\mathcal{F})$ a \emph{Dirichlet form} on $L^2(M,\mu)$, if
\begin{itemize}
	\item $\mathcal{F}$ is a Hilbert space equipped with the inner product $\mathcal{E}_1(u,v):=\mathcal{E}(u,v)+(u,v)_{L^2(\mu)}$;
	\item the Markovian property holds: $f_+\wedge 1\in\mathcal{F}$ for every $f\in\mathcal{F}$, and $\mathcal{E}(f_+\wedge 1,f_+\wedge 1)\le\mathcal{E}(f,f)$.
\end{itemize}
It is called \emph{regular}, if $\mathcal{F}\cap C_0(M,d)$ is dense both in $\mathcal{F}$ (under norm $\mathcal{E}_1^{1/2}$) and in $C_0(M,d)$ (under norm $\|\cdot\|_\infty$). In this case, for any open set $\Omega\subset M$ and Borel subset $U\Subset\Omega$, there exists $\phi\in\mathcal{F}\cap C_0(\Omega)$ such that $\phi=1$ on $U$ and $0\le\phi\le 1$ throughout $M$. The collection of such $\phi$ is denoted by $\mathrm{cutoff}(U,\Omega)$.

Recall that every regular Dirichlet form admits the following Beurling-Deny decomposition:
\begin{equation}\label{bd}
	\mathcal{E}(u,v)=\mathcal{E}^{(L)}(u,v)+\int_{M\times M}(u(x)-u(y))(v(x)-v(y))dj(x,y)+\int_Mu(x)v(x)d\kappa(x)
\end{equation}
for all $u,v\in\mathcal{F}\cap C_0(\Omega)$, where $\mathcal{E}^{(L)}$ is \emph{strongly local}, that is, $\mathcal{E}^{(L)}(u,v)=0$ whenever $\mathrm{supp}(u-a)\cap\mathrm{supp}(v)=\emptyset$ with some constant $a\in\mathbb{R}$; the \emph{jump measure} $j$ is a symmetric Radon measure on $M\times M$ with $j(\mathrm{diag})=0$; the \emph{killing measure} $\kappa$ is a Radon measure on $M$. See \cite[Theorem 4.5.2]{fot} for the meaning of these notions.

Throughout this paper we assume that $\kappa\equiv 0$, and
$$j(dx,dy)=d\mu(x)J(x,dy)=d\mu(y)J(y,dx)$$
with a Borel transition kernel $J$. In this case, $\mathcal{E}$ can be extended to $\mathcal{F}'$, the collection of all functions $u=v+a$, where $v\in\mathcal{F}$ and $a\in\mathbb{R}$, by $\mathcal{E}(u,u):=\mathcal{E}(v,v)$. Further, (\ref{bd}) holds for all $u,v\in\mathcal{F}'$.

Recall that for any $u\in\mathcal{F}$, there is a unique Radon measure $\Gamma(u)$ on $M$, called the (le Jan) \emph{energy form} of $u$, determined by
$$\int_M\varphi d\Gamma(u)=\mathcal{E}(u,\varphi u)-\frac{1}{2}\mathcal{E}(u^2,\varphi)\quad\mbox{for all}\quad\varphi\in\mathcal{F}.$$
Under our assumptions, $\Gamma(u)$ can also be defined for any $u\in\mathcal{F}'$ (with $\Gamma(c)\equiv 0$ for any constant $c$), and
$$d\Gamma(u)(x)=d\Gamma^{(L)}(u)(x)+\left(\int_M(u(x)-u(y))^2J(x,dy)\right)d\mu(x),$$
where $\Gamma^{(L)}$ is the energy form of $u$ with respect to $\mathcal{E}^{(L)}$. It defines a symmetric bilinear form $\Gamma^{(L)}(u,v)$ which satisfies the product and Leibniz rules by strong locality (cf.\ \cite[Lemma 3.2.5 and Theorem 3.2.2]{fot}).

For every open set $\Omega\subset M$, we also define the $\Omega$-local energy
$$d\Gamma_\Omega(u)(x)=1_\Omega(x)\left\{d\Gamma^{(L)}(u)(x)+\left(\int_\Omega(u(x)-u(y))^2J(x,dy)\right)d\mu(x)\right\}.$$
Let $\mathcal{F}(\Omega)$ be the closure of $\mathcal{F}\cap C_0(\Omega)$ in $(\mathcal{F},\mathcal{E}_1)$, and $P_t^\Omega$ be the \emph{heat semigroup} on $L^2(\Omega)$ associated with the Dirichlet form $(\mathcal{E},\mathcal{F})$, that is, for any $f\in L^2(\Omega)$ and $\varphi\in\mathcal{F}(\Omega)$,
$$\lim\limits_{t\downarrow 0}\left\|P_t^\Omega f-f\right\|_{L^2(\Omega)}=0,\quad\mbox{and}\quad\frac{d}{dt}\left(P_t^\Omega f,\varphi\right)+\mathcal{E}\left(P_t^\Omega f,\varphi\right)=0\quad\mbox{for all}\quad t>0.$$
It can always be extended to a contractive semigroup on $L^p(\Omega)$ with $1\le p\le\infty$, still denoted as $P_t^\Omega$. Clearly there exists a transition kernel $p_t^\Omega(x,dy)$ on $\Omega$ such that for $\mu$-a.e.\ $x\in\Omega$, all $f\in L^p(\Omega)$ with any $1\le p\le\infty$,
$$P_t^\Omega f(x)=\int_\Omega f(y)p_t^\Omega(x,dy).$$
Be careful that $\Gamma_\Omega$ is the energy form of the Neumann boundary problem on $\Omega$, but $P_t^\Omega$ is the heat semigroup of the Dirichlet boundary problem on $\Omega$. Generally they do not correspond to each other.

If $p_t^\Omega(x,dy)\ll d\mu(y)$, then we call the Radon-Nykodym derivative $p_t^\Omega(x,y)$ as the \emph{heat kernel} of $P_t^\Omega$. When $\Omega=M$, we simply write $P_t$ and $p_t(x,y)$ (if exists) instead of $P_t^M$ and $p_t^M(x,y)$, respectively called the (global) heat semigroup and heat kernel associated with the Dirichlet form $(\mathcal{E},\mathcal{F})$.

Let $\mathcal{F}'(\Omega)$ be the collection of all $v\in L_{loc}^1(M)$ such that there exist $w\in\mathcal{F}$ and a strict neighbourhood $\Omega'$ of $\Omega$ (that is, $U\cap\Omega\Subset\Omega'$ for any bounded open set $U\subset M$) such that $v=w$ on $\Omega'$ and
$$\int_\Omega\int_M(v(x)-v(y))^2J(x,dy)d\mu(x)<\infty.$$
For all $v\in\mathcal{F}'(\Omega)$ and $\varphi\in\mathcal{F}(\Omega)$, let the extended Dirichlet form be
$$\tilde{\mathcal{E}}(v,\varphi):=\mathcal{E}^{(L)}(w,\varphi)+\int_{M\times M}(\varphi(x)-\varphi(y))(v(x)-v(y))dj(x,y).$$

Let $Q=(t_1,t_2]\times\Omega$ be an arbitrary cylinder in $\mathbb{R}\times M$, where $t_1<t_2$, and $\Omega$ is an open subset of $M$. For any function $u:(t_1,t_2]\to\mathcal{F}'(\Omega)$, let
\begin{equation}\label{e'sup}
\mathop{\mathrm{esup}}_Qu:=\sup\limits_{t_1<t\le t_2}\mathop{\mathrm{esup}}_\Omega u(t,\cdot),\quad\mathop{\mathrm{einf}}_Qu:=\inf\limits_{t_1<t\le t_2}\mathop{\mathrm{einf}}_\Omega u(t,\cdot).
\end{equation}

On $\mathbb{R}\times M$, let $\hat{\mu}$ be the completion of the product measure $m_1\otimes\mu$. Given two open sets $U,\Omega$ in $M$ with $U\Subset\Omega$, define the tail of a $j$-measurable function $f$ with respect to $(U,\Omega)$ as
$$T_\Omega^U(f)=\mathop{\mathrm{esup}}\limits_{z\in U}\int_{\Omega^c}|f(y)|J(z,dy).$$
Set $Q=(t_1,t_2]\times\Omega$ with some $t_1<t_2$. For any $\hat{\mu}$-measurable function $u$ on $(t_1,t_2]\times M$, we define $T_Q^{U;p}(u)$, the $L^p$-tail of $u$ with respect to $(U,Q)$, as the $L^p$-norm of $t\mapsto T_\Omega^U(u(t,\cdot))$ on $(t_1,t_2]$, that is,
$$T_Q^{U;p}(u)=\left(\int_{t_1}^{t_2}\left(T_\Omega^U(u(s,\cdot))\right)^pds\right)^{1/p}\quad\mbox{for}\quad 1\le p<\infty,\quad\mbox{and}\quad T_Q^{U;\infty}(u)=\mathop{\mathrm{esup}}\limits_{t_1<s\le t_2}T_\Omega^U(u(s,\cdot)).$$
If $\Omega=B$ and $U=\lambda B$ with some metric ball $B$ and $0<\lambda<1$, then we also write these tails as $T_B^{(\lambda)}(f)$ and $T_Q^{(\lambda;p)}(u)$ for simplicity.

We say $u:(t_1,t_2]\to\mathcal{F}'(\Omega)$ is \emph{caloric} on $Q$, if
\begin{itemize}
\item $u$ is $\hat{\mu}$-measurable on $(t_1,t_2]\times M$;
\item for any $\varphi\in L^2(\Omega)$, $I_\varphi(t):=(u(t,\cdot),\varphi)_{L^2(\Omega)}$ is continuous on $(t_1,t_2]$, and the left derivative $\frac{d_-}{dt}I_\varphi(t)$ exists at every $t_1<t\le t_2$. Further, if $\varphi\in\mathcal{F}(\Omega)$, then
\begin{equation}\label{cal}
\frac{d_-}{dt}\big(u(t,\cdot),\varphi\big)+\tilde{\mathcal{E}}(u(t,\cdot),\varphi)=0;
\end{equation}
\item there exist $t_1<s_0<\cdots<s_N=t_2$ (where $N<\infty$) such that $\frac{d}{dt}I_\varphi(t)$ exists for all $t\in(t_1,t_2]\setminus\{s_1,\dotsc,s_N\}$ and $\varphi\in L^2(\Omega)$.
\end{itemize}

Same as we introduced in \cite{L0}, this definition covers \cite[Definition 2.8]{ghh+2}, and a caloric function $u$ satisfies
\begin{equation}\label{cal'}
\big(u(t,\cdot),\varphi\big)=\lim_{s\downarrow t_1}\big(u(s,\cdot),\varphi\big)-\int_{t_1}^t\tilde{\mathcal{E}}(u(s,\cdot),\varphi)ds
\end{equation}
for all $\varphi\in\mathcal{F}(\Omega)$ and $t_1<t\le t_2$.

Let $W:M\times\mathbb{R}_+\to\mathbb{R}_+$ be a \emph{scaling function} (for short, a \emph{scale}), that is, given any $x\in M$, $W(x,\cdot)$ is strictly increasing, continuous, and there exist $C_W\ge 1$ and $0<\beta_1\le\beta_2<\infty$ such that for all $0<r\le R<\infty$ and $x,y\in M$ with $d(x,y)\le R$,
\begin{equation}\label{W}
C_W^{-1}\left(\frac{R}{r}\right)^{\beta_1}\le\frac{W(x,R)}{W(y,r)}\le C_W\left(\frac{R}{r}\right)^{\beta_2}.
\end{equation}
It follows clearly that $W(x,0)=0$ and $w(x,R)\to\infty$ as $R\to\infty$ for any $x\in M$.

For a ball $B=B(x,r)$, we also denote $W(B):=W(x,r)$. For any $t\ge 0$, let $W^{-1}(x,t)$ be the unique $r\in\mathbb{R}_+$ such that $t=W(x,r)$. It follows obviously by (\ref{W}) that for all $0<s\le t<\infty$ and $x\in M$,
\begin{equation}\label{W-1}
\left(\frac{t}{C_Ws}\right)^{\frac{1}{\beta_2}}\le\frac{W^{-1}(x,t)}{W^{-1}(x,s)}\le\left(\frac{C_Wt}{s}\right)^{\frac{1}{\beta_1}}.
\end{equation}

Fixing a scale $W$, we define the following conditions (some other conditions are omitted but will be mentioned in Corollary \ref{1-18} below. They have been completely investigated in \cite{L0}. See the definitions there):

\begin{itemize}
\item $(\mathrm{PHI}^0)$, the basic parabolic Harnack inequality: there exist constants $0<\delta_1<\delta_2<\delta_3<\delta_4<\infty$, $C\ge 1$, $\lambda_0\in(0,1)$ and $\lambda_1>0$ such that for every ball $B=B(x_0,R)$ with $0<R<\overline{R}$ and all $u:(0,\delta_4W(\lambda_1B)]\to\mathcal{F}'(B)$ that is non-negative on $(0,\delta_4W(\lambda_1B)]\times M$, bounded and caloric on $Q=(0,\delta_4W(\lambda_1B)]\times B$,
\begin{equation}\label{phi0}
	\mathop{\mathrm{esup}}_{(\delta_1W(\lambda_1B),\delta_2W(\lambda_1B)]\times\lambda_0B}u\le C\mathop{\mathrm{einf}}_{(\delta_3W(\lambda_1B),\delta_4W(\lambda_1B)]\times\lambda_0B}u.
\end{equation}

\item $(\mathrm{PHI})$, the standard parabolic Harnack inequality: there exist $\delta>0$, $C\ge 1$ and $\lambda_0\in(0,1)$ such that for every ball $B=B(x_0,R)$ with $0<R<\overline{R}$ and all $u:(t_0-4\delta W(B),t_0]\to\mathcal{F}'(B)$ that is non-negative, bounded and caloric on $Q=(t_0-4\delta W(B),t_0]\times B$,
\begin{equation}\label{phi}
	\mathop{\mathrm{esup}}_{(t_0-3\delta W(B),t_0-2\delta W(B)]\times\lambda_0B}u\le C\left(\mathop{\mathrm{einf}}_{(t_0-\delta W(B),t_0]\times\lambda_0B}u+T_Q^{(3/4;1)}(u_-)\right).
\end{equation}

\item $(\mathrm{PHI}^+)$, the complete parabolic Harnack inequality: there exist $\delta>0$, $C\ge 1$ and $\lambda_0\in(0,1)$ such that for every ball $B=B(x_0,R)$ with $0<R<\overline{R}$ and all $u:(t_0-4\delta W(B),t_0]\to\mathcal{F}'(B)$ that is non-negative, bounded and caloric on $Q=(t_0-4\delta W(B),t_0]\times B$, (with $Q_\pm$ defined same as in $(\mathrm{PHI})$ condition)
\begin{equation}\label{phi+}
	\mathop{\mathrm{esup}}_{Q_-}u+T_{Q_-}^{(2/3;1)}(u_+)\le C\left(\mathop{\mathrm{einf}}_{Q_+}u+T_Q^{(3/4;1)}(u_-)\right),
\end{equation}
where $Q_-:=(t_0-3\delta W(B),t_0-2\delta W(B)]\times\lambda_0B$ and $Q_+:=(t_0-\delta W(B),t_0]\times\lambda_0B$.

Such an inequality appeared first in \cite[Corollary 1.4]{kw} for classical $\alpha$-stable processes.

\item $(\mathrm{LLE})$, the local lower estimate (of the heat kernel) (also denoted as $(\mathrm{NDL})$ in other literature): there exist $\delta_L,c_L\in(0,1)$ such that for every ball $B=B(x_0,R)$ with $0<R<\overline{R}$, the heat kernel $p_t^B(x,y)$ exists, and for all $0<t\le W(\delta_LB)$,
$$p_t^B(x,y)\ge\frac{c_L}{V\left(x_0,W^{-1}(x_0,t)\right)}\quad\quad\mu\mbox{-a.e.}\quad\quad x,y\in B\left(x_0,\delta_LW^{-1}(x_0,t)\right).$$

\item $(\mathrm{UE})$, the upper estimate (of the heat kernel): the heat kernel $p_t(x,y)$ exists, and there exists $C>0$ such that for $\mu$-a.e.\ $x,y\in M$ and all $0<t<W(x,\overline{R})$,
$$p_t(x,y)\le C\left(\frac{1}{V\left(x,W^{-1}(x,t)\right)}\wedge\frac{t}{V(x,d(x,y))W(x,d(x,y))}\right).$$

\item $(\mathrm{PMV}^-)$, the parabolic mean-value inequality: there exists a constant $C>0$ such that for all $B=B(x_0,R)$ with $0<R<\overline{R}$ and all $\varepsilon>0$, if $u:(t_0-W(B),t_0]\to\mathcal{F}'(B)$ is non-negative, bounded and subcaloric on $Q=(t_0-W(B),t_0]\times B$, then
\begin{equation}\label{pmv2p}
\mathop{\mathrm{esup}}\limits_{Q_+}u\le C\left(1+\varepsilon^{-\frac{1+\nu}{2\nu}}\right)\left(\fint_{Q'}u^2d\hat{\mu}\right)^{1/2}+\varepsilon\sup\limits_{t_0-\frac{1}{2}W(B)\le s\le t_0}\|u_+(s,\cdot)\|_{L^\infty(M\setminus\frac{1}{2}B)},
\end{equation}
where $Q':=(t_0-\frac{1}{2}W(B),t_0]\times B$ and $Q_+:=(t_0-\frac{1}{4}W(B),t_0]\times\frac{1}{2}B$.

This condition was referred to as $(\mathrm{PMV}_1)$ in \cite{ghh+2}.

\item $(\mathrm{PMV}^+)$, the strengthened parabolic mean-value inequality: there exist $C,\delta>0$ and $c_0\in(0,1)$ such that for all $B=B(x_0,R)$ with $0<R<\overline{R}$, if $u:(t_0-\delta W(B),t_0+\frac{1}{4}\delta W(B)]\to\mathcal{F}'(B)$ is bounded, caloric on $Q=(t_0-\delta W(B),t_0+\frac{1}{4}\delta W(B)]\times B$ and non-negative on $(t_0-\delta W(B),t_0+\frac{1}{4}\delta W(B)]\times M$, then
\begin{equation}\label{p1mv-p}
\mathop{\mathrm{esup}}\limits_{(t_0-\frac{1}{4}\delta W(B),t_0]\times c_0B}u\le C\fint_{(t_0-\frac{1}{2}\delta W(B),t_0+\frac{1}{4}\delta W(B)]\times B}ud\hat{\mu}.
\end{equation}

\item $(\mathrm{wPH}_1)$, the $L^1$-weak parabolic Harnack inequality: there exist $\delta_0>0$, $C>1$ and $\lambda_0\in(0,1)$ such that for every ball $B=B(x_0,R)$ with $R\in(0,\overline{R})$, all $0<\lambda\le\lambda_0$ and all $u$ that is non-negative, bounded and supercaloric on $Q=(t_0-4\delta_0W(\lambda B),t_0]\times B$,
\begin{equation}\label{wph'q}
\fint_{(t_0-3\delta_0W(\lambda B),t_0-2\delta_0W(\lambda B)]\times\lambda B}ud\hat{\mu}\le C\left(\mathop{\mathrm{einf}}_{(t_0-\delta_0W(\lambda B),t_0]\times\lambda B}u+T_Q^{(2/3;1)}(u_-)\right).
\end{equation}

\item $(\mathrm{PGL}_1)$, the $L^1$-parabolic growth lemma: there exist $\delta_0>0$, $C>1$ and $\lambda_0\in(0,1)$ such that for every ball $B=B(x_0,R)$ with $R\in(0,\overline{R})$, all $0<\lambda\le\lambda_0$ and all $u$ that is non-negative, bounded and supercaloric on $Q=(t_0-4\delta_0W(\lambda B),t_0]\times B$, if
$$\frac{\hat{\mu}\big(Q_-\cap\{u\ge a\}\big)}{\hat{\mu}(Q_-)}\ge\eta$$
with some $a,\eta>0$, where $Q_-=(t_0-3\delta_0W(\lambda B),t_0-2\delta_0W(\lambda B)]\times\lambda B$, then
$$\mathop{\mathrm{einf}}_{(t_0-\delta_0W(\lambda B),t_0]\times\lambda B}u\ge C^{-1}\eta a-2T_Q^{(2/3;1)}(u_-).$$

\item $(\mathrm{cap}_\le)$, the upper bound of capacity: for any $0<\lambda<1$, there exists $C_\lambda>0$ such that for every ball $B=B(x,R)$ with $0<R<\overline{R}$, there exists $\phi\in\mathrm{cutoff}(\lambda B,B)$ such that
$$\mathcal{E}(\phi,\phi)\le C_\lambda\frac{\mu(B)}{W(B)}.$$

\item $(\mathrm{UJS})$, the upper jumping smoothness: the jump kernel $J:M\times M\setminus\mathrm{diag}\to\mathbb{R}_+$ exists (that is, the measure $J(x,dy)=J(x,y)d\mu(y)$ for $\mu$-a.e. $x\in M$), and there exists $C>0$ such that for $\mu$-a.e. $x,y\in M$,
$$J(x,y)\le\frac{C}{V(x,r)}\int_{B(x,r)}J(z,y)d\mu(z)\quad\mbox{with all}\quad 0<r\le\frac{1}{2}d(x,y).$$

\item $(\mathrm{TJ})$, the tail of jump measure: for any $0<\lambda<1$, there exists $C_\lambda>0$ such that for every ball $B=B(x_0,R)$,
$$\mathop{\mathrm{esup}}\limits_{x\in\lambda B}\int_{B^c}J(x,dy)\le\frac{C_\lambda}{W(B)}.$$
\end{itemize}
For other related conditions that we will mention literally below, see \cite[Section 2]{L0} for definitions. Recall that
\begin{equation}\label{3-1}
(\mathrm{TJ})+(\mathrm{UJS})\Rightarrow(\mathrm{J}_\le)\quad\mbox{and}\quad(\mathrm{J}_\le)+(\mathrm{VD})\Rightarrow(\mathrm{TJ}),
\end{equation}
while by \cite[Theorem 2.1]{L0},
\begin{equation}\label{M1}
(\mathrm{LLE})\Leftrightarrow(\mathrm{wPH}_1)\Rightarrow(\mathrm{PGL}_1).
\end{equation}

The target of this paper is to prove the following main theorem (which we hint in \cite[Corollary 2.6]{L0} using probability) by pure analysis:

\begin{theorem}\label{M4}
	Assume that $(M,d,\mu)$ satisfies $(\mathrm{VD})$ and $(\mathrm{RVD})$, and $(\mathcal{E},\mathcal{F})$ is a regular Dirichlet form satisfying (\ref{bd}). Then
	$$(\mathrm{wPH}_1)+(\mathrm{UJS})\Leftrightarrow(\mathrm{PHI}^0)\Leftrightarrow(\mathrm{PHI})\Leftrightarrow(\mathrm{PHI}^+).$$
\end{theorem}

\begin{proof}
We show $(\mathrm{wPH}_1)+(\mathrm{UJS})\Rightarrow(\mathrm{PHI}^0)$ by Theorem \ref{main''}; $(\mathrm{PHI}^0)\Rightarrow(\mathrm{PHI})\Rightarrow(\mathrm{PHI}^+)$ by Propositions \ref{phi+t} and \ref{ph+}; $(\mathrm{PHI}^+)\Rightarrow(\mathrm{wPH}_1)+(\mathrm{UJS})$ by a trivial implication and Proposition \ref{ujs}.
\end{proof}

Combining all the aforementioned ideas, we have the following equivalent characterizations of $(\mathrm{PHI})$, greatly enlarging the list of equivalent conditions raised by \cite{ckw-e,ckw-p}:

\begin{corollary}\label{1-18}
	Under the same assumptions as in Theorem \ref{M4}, each of the following groups of conditions is equivalent to $(\mathrm{PHI})$:
	\begin{longtable}{clccl}
		$\mathrm{(a)}$&$(\mathrm{PHI}^+)$&\quad\quad\quad&$\mathrm{(a)'}$&$(\mathrm{PHI}^0)$\\
		$\mathrm{(b)}$&$(\mathrm{EHI})+(\mathrm{E})+(\mathrm{UJS})$&&$\mathrm{(b)'}$&$(\mathrm{EHI}^0)+(\mathrm{E})+(\mathrm{UJS})$\\
		$\mathrm{(c)}$&$(\mathrm{wPH}_1)+(\mathrm{UJS})$&&$\mathrm{(c)'}$&$(\mathrm{wPH})+(\mathrm{FK})+(\mathrm{UJS})$\\
		$\mathrm{(d)}$&$(\mathrm{wPH}_1)+(\mathrm{PMV}^+)$&&$\mathrm{(d)'}$&$(\mathrm{wPH})+(\mathrm{PMV}^+)$\\
		$\mathrm{(e)}$&$(\mathrm{wPH}_1)+(\mathrm{PMV}^-)+(\mathrm{UJS})$&&$\mathrm{(e)'}$&$(\mathrm{PMV}^-)+(\mathrm{cap}_\le)+(\mathrm{PI})+(\mathrm{UJS})$\\
		$\mathrm{(f)}$&$(\mathrm{wEH})+(\mathrm{E})+(\mathrm{UJS})$&&$\mathrm{(f)'}$&$(\mathrm{wEH})+(\mathrm{cap}_\le)+(\mathrm{FK})+(\mathrm{UJS})$\\
		$\mathrm{(g)}$&$(\mathrm{PHR})+(\mathrm{E})+(\mathrm{UJS})$&&$\mathrm{(g)'}$&$(\mathrm{PHR}^0)+(\mathrm{E})+(\mathrm{UJS})$\\
		$\mathrm{(h)}$&$(\mathrm{EHR})+(\mathrm{E})+(\mathrm{UJS})$&&$\mathrm{(h)'}$&$(\mathrm{EHR}^0)+(\mathrm{E})+(\mathrm{UJS})$\\
		$\mathrm{(j)}$&$(\mathrm{LLE})+(\mathrm{UE})+(\mathrm{UJS})$&&$\mathrm{(j)'}$&$(\mathrm{LLE})+(\mathrm{UJS})$\\
		$\mathrm{(k)}$&$(\mathrm{PI})+(\mathrm{E})+(\mathrm{UJS})$&&$\mathrm{(k)'}$&$(\mathrm{PI})+(\mathrm{S}^-)+(\mathrm{UJS})$\\
		$\mathrm{(l)}$&$(\mathrm{PI})+(\mathrm{Gcap})+(\mathrm{J}_\le)+(\mathrm{UJS})$&&$\mathrm{(l)'}$&$(\mathrm{PI})+(\mathrm{Gcap})+(\mathrm{UJS})$\\
		$\mathrm{(m)}$&$(\mathrm{PI})+(\mathrm{ABB})+(\mathrm{J}_\le)+(\mathrm{UJS})$&&$\mathrm{(m)'}$&$(\mathrm{PI})+(\mathrm{ABB})+(\mathrm{TJ})+(\mathrm{UJS})$
	\end{longtable}
\end{corollary}

Note that the implications $\mathrm{(\bullet)}\Rightarrow\mathrm{(\bullet)'}$ are trivially true.

\begin{proof}
	$(\mathrm{PHI})\Leftrightarrow\mathrm{(a)}\Leftrightarrow\mathrm{(a)'}\Leftrightarrow\mathrm{(c)}$ by Theorem \ref{M4}. Further, $\mathrm{(c)}\Leftrightarrow\mathrm{(d)}\Leftrightarrow\mathrm{(e)}$ by Theorem \ref{main''} and its proof.
	
	By (\ref{M1}) we see $\mathrm{(c)}\Leftrightarrow\mathrm{(c)'}\Leftrightarrow\mathrm{(j)'}\Rightarrow(\mathrm{cap}_\le)$, which implies $(\mathrm{TJ})$ by Lemma \ref{L2-J} since $(\mathrm{UJS})$ is contained here. Thus $\mathrm{(j)'}\Leftrightarrow\mathrm{(f)}\Leftrightarrow\mathrm{(f)'}\Leftrightarrow\mathrm{(g)}\Leftrightarrow\mathrm{(g)'}\Leftrightarrow\mathrm{(h)}\Leftrightarrow\mathrm{(h)'}\Leftrightarrow\mathrm{(k)}\Leftrightarrow\mathrm{(k)'}\Leftrightarrow\mathrm{(m)'}$ by \cite[Theorem 2.2]{L0}
	
	Trivially $\mathrm{(d)}\Rightarrow\mathrm{(d)'}$; $\mathrm{(d)'}\Rightarrow\mathrm{(a)'}$ by (\ref{pmv'nas}); $\mathrm{(e)}\Rightarrow\mathrm{(e)'}$ by \cite[Theorem 2.2]{L0} and the trivial implication $(\mathrm{Gcap})\Rightarrow(\mathrm{cap}_\le)$; $\mathrm{(e)'}\Leftrightarrow\mathrm{(k)}$ by \cite[Corollaries 10.3, 12.5 and Lemma 11.2]{ghh+1}. Hence $\mathrm{(a)'}\Leftrightarrow\mathrm{(d)'}\Leftrightarrow\mathrm{(e)'}$.
	
	Trivially $\mathrm{(a)'}+\mathrm{(h)'}\Rightarrow\mathrm{(b)'}$, while $\mathrm{(b)'}\Rightarrow\mathrm{(b)}$ by \cite[Proposition 5.4]{L0}; $\mathrm{(b)}\Rightarrow\mathrm{(h)}$ in the same way as in \cite[Lemma 6.2]{hy}. Thus we also have $\mathrm{(h)'}\Leftrightarrow\mathrm{(b)}\Leftrightarrow\mathrm{(b)'}$.
	
	Finally, $\mathrm{(m)'}\Rightarrow\mathrm{(m)}$ by Lemma \ref{L2-J} (the second assertion); $\mathrm{(m)}\Rightarrow\mathrm{(l)'}$ by the implication $(\mathrm{iii})'\Rightarrow(\mathrm{iii})''$ in \cite[Theorem 2.2]{L0} (since $(\mathrm{TJ})$ holds by (\ref{3-1})); $\mathrm{(l)'}\Rightarrow\mathrm{(l)}$ by Lemma \ref{L2-J}; $\mathrm{(l)}\Rightarrow\mathrm{(j)}$ by \cite[Subsection 8.2]{ghh+3}, and $\mathrm{(j)}\Rightarrow\mathrm{(j)'}$ trivially. Therefore, $\mathrm{(m)'}\Leftrightarrow\mathrm{(m)}\Leftrightarrow\mathrm{(l)}\Leftrightarrow\mathrm{(l)'}$. The proof is then completed.
\end{proof}

In the sequel, $C$ always stands for a positive constant (independent of the variables in the context but possibly depending on parameters), whose value may change from line to line, but is typically greater than 1. When necessary, we refer to the constant $C$ coming from condition $(\bullet)$ by $C_\bullet$ so that its value is fixed once appearing.

\section{Basic properties}
\subsection{Comparison inequalities}

\begin{lemma}\cite[Proposition 3.3]{L0}\label{gpmp}
	Let $\Omega$ be a bounded open set in $M$, and $Q=(t_1,t_2]\times\Omega$. If
	\begin{itemize}
		\item $u:(t_1,t_2]\to\mathcal{F}'(\Omega)$ is subcaloric on $Q$;
		\item $u\le 0$ weakly on $(t_1,t_2]\times\Omega^c$, that is, $u_+(t,\cdot)\in\mathcal{F}(\Omega)$ for all $t_1<t\le t_2$;
		\item $u|_{t=t_1}\le 0$ weakly in $\Omega$, that is, $\|u_+(t,\cdot)\|_{L^2(\Omega)}\to 0$ as $t\downarrow t_1$,
	\end{itemize}
	then $u\le 0$ on $(t_1,t_2]\times M$.
\end{lemma}

It is easy to check that:

\begin{corollary}\label{aux-comp}
	For any $u$ that is supercaloric on $Q=(t_1,t_2]\times\Omega$ and non-negative on $(t_1,t_2]\times M$, for any measurable $K$ and open $U$ such that $K\Subset U\Subset\Omega$, any $t_1<s<t\le t_2$ and $\mu$-a.e.\ $x\in M$,
	$$u(t,x)\ge\mathop{\mathrm{einf}}_{z\in K}u(s,z)\cdot P_{t-s}^U1_K(x).$$
\end{corollary}

The following lemma is extremely important throughout this paper:

\begin{lemma}\label{trunc}
	Given any open sets $U\Subset\Omega$ and any $\phi\in\mathrm{cutoff}(U,\Omega)$, for every $u:(t_1,t_2]\to\mathcal{F}'(\Omega)$ that is subcaloric on $Q=(t_1,t_2]\times\Omega$ with $T_Q^{U;1}(u_+)<\infty$, then
	$$u'_{U;\Omega}(t,x):=u(t,x)\phi(x)+2T_{(t,t_2]\times\Omega}^{U;1}(u_+)$$
	is subcaloric on $Q'=(t_1,t_2]\times U$.
\end{lemma}

\begin{proof}
	For every non-negative $\varphi\in\mathcal{F}(U)$, we have
	\begin{align*}
		\frac{d}{dt}\left(u'_{U;\Omega},\varphi\right)=&\ \frac{d}{dt}(u\phi,\varphi)-2\left(T_\Omega^U(u_+(t,\cdot)),\varphi\right)=\frac{d}{dt}(u\phi,\varphi)-2\left(\mathop{\mathrm{esup}}\limits_{x\in U}\int_{\Omega^c}u_+(t,z)J(x,dz)\right)\int_U\varphi d\mu\\
		\le&\ \frac{d}{dt}(u,\varphi)-2\int_{U\times\Omega^c}\varphi(x)u_+(t,z)dj(x,z)\le-\mathcal{E}(u,\varphi)-2\int_{U\times\Omega^c}\varphi(x)u(t,z)(1-\phi_n(z))dj(x,z)\\
		=&\ -\mathcal{E}^{(L)}(u,\varphi)+2\int_{U\times\Omega^c}\varphi(x)u(t,z)\phi(z)dj(x,z)=-\mathcal{E}(u\phi,\varphi)=-\mathcal{E}\left(u'_{U;\Omega},\varphi\right)
	\end{align*}
	for all $t_1<t\le t_2$, which implies that $u'_{U;\Omega}$ is subcaloric on $Q'$.
\end{proof}

\subsection{Semi-homogeneous caloric exension}

Given every two open sets $U,\Omega$ in $M$ with $\Omega\Supset U$, if $f:(t_1,t_2]\to L^1_{loc}(M,\mu)$ is $\hat{\mu}$-measurable and non-negative on $(t_1,t_2]\times M$, and $\mathrm{supp}(f(t,\cdot))\subset\Omega^c$ for all $t_1<t\le t_2$, then we could define the \emph{caloric extension} of $f$ into $(t_1,t_2]\times U$ as
\begin{equation}\label{ext-d}
f_{t_1,U}(t,x)=\begin{cases}f(t,x),&\mbox{if }x\in U^c;\\
2\int_{t_1}^t\int_U\int_{\Omega^c}f(s,y)J(z,dy)p_{t-s}^U(x,dz)ds,&\mbox{if }x\in U
\end{cases}\end{equation}
on $(t_1,t_2]\times M$, which is also valid for any real-valued $f$ provided the integration is well-defined.

\begin{lemma}\label{ext-F}\cite[Lemma 5.2]{L0}
For every ball $B=B(x_0,R)$ and every $0<\lambda<1$, if $f$ is $\hat{\mu}$-measurable, non-negative and $\mathrm{supp}f(t,\cdot)\subset B^c$ for all $t_1<t\le t_2$,

(i) if $T_{(t_1,t_2]\times B}^{(\lambda;1)}(f)<\infty$, then $1_Bf_{t_1,\lambda B}$ is bounded on $(t_1,t_2]\times\lambda B$, and for $\mu$-a.e. $x\in B$,
\begin{equation}\label{ext-infty}
f_{t_1,\lambda B}(t,x)\le 2T_{(t_1,t]\times B}^{(\lambda;1)}(f);
\end{equation}

(ii) if $T_{(t_1,t_2]\times B}^{(\lambda;2)}(f)<\infty$, then $1_Bf_{t_1,\lambda B}(t,\cdot)\in\mathcal{F}(\lambda B)$ for all $t_1<t\le t_2$;

(iii) if further $f(t,\cdot)\in\mathcal{F}'(\lambda B)$ for all $t_1<t\le t_2$, then $f_{t_1,\lambda B}\in\mathcal{F}'(\lambda B)$ is caloric on $(t_1,t_2]\times\lambda B$;

(iv) if $f(t,\cdot)\in\mathcal{F}'$ for all $t$, then $f_{t_1,\lambda B}\in\mathcal{F}'$; while if $f(t,\cdot)\in\mathcal{F}$ for all $t$, then $f_{t_1,\lambda B}\in\mathcal{F}$.
\end{lemma}

In Subsection \ref{p'ujs} we also need the caloricity of $f_{t_1,U}$ in a different setting, for which we need some preparation:

\begin{lemma}
	For any open sets $U\Subset\Omega$ in $M$, any $0<t<\infty$ and $\mu$-a.e.\ $z\in U$,
	\begin{equation}\label{levy+}
		2\int_0^t\int_{U\times\overline{\Omega}^c}p_s^\Omega(x,z)dj(x,y)ds\le\left(1-P_t^\Omega1_\Omega(z)\right)\le 1.
	\end{equation}
\end{lemma}

In probability theory, the same could be proved by the L\'evy system formula, see the first three lines in \cite[Proof of Corollary 3.4]{ckw-p} -- the coefficient ``2'' disappears there simply because the definition of the jump measure is different (see \cite[Appendix B]{L0}).

\begin{proof}
	For any compact $K\subset\overline{\Omega}^c$, let $\Phi\in\mathrm{cutoff}\left(K,\overline{\Omega}^c\right)$. By \cite[Lemma 4.1]{ghh17},
	\begin{align*}
		\left(1-P_t^\Omega1_\Omega,f\right)\ge&\ \int_0^t-\mathcal{E}(\Phi,P_s^\Omega f)ds=\int_0^t\left\{2\int_{\Omega\times\overline{\Omega}^c}P_s^\Omega f(x)\Phi(y)dj(x,y)\right\}ds\\
		=&\ 2\int_0^t\int_\Omega\int_{\Omega\times\overline{\Omega}^c}p_s^\Omega(x,z)\Phi(y)dj(x,y)f(z)d\mu(z)ds\ge 2\int_0^t\int_U\int_{U\times K}p_s^\Omega(x,z)dj(x,y)f(z)d\mu(z)ds
	\end{align*}
	for every non-negative $f\in L^1\cap L^2$ supported in $U$. Letting $K$ exhaust $\overline{\Omega}^c$, we have
	$$\int_U\left(\int_0^t\int_{U\times\overline{\Omega}^c}p_s^\Omega(x,z)dj(x,y)ds-\frac{1}{2}\left(1-P_t^\Omega1_\Omega(z)\right)\right)f(z)d\mu(z)\le 0.$$
	The desired inequality follows since $f$ is arbitrary.
\end{proof}

\begin{proposition}\label{cal+}
	Let $\Omega$ be a precompact open subsets of $M$ and $K$ be a compact set with $K\cap\overline{\Omega}=\emptyset$. Then for any non-negative $f\in L^\infty$ supported in $K$ and any $t_1\le t_3<t_4<\infty$, $(f1_{(t_3,t_4]})_{t_1,\Omega}$ is caloric on $(t_1,\infty)\times\Omega$.
\end{proposition}

\begin{proof}
	Fix $\phi\in\mathrm{cutoff}(K,\overline{\Omega}^c)$. For simplicity, we denote
	$$T_f(z)=\int_Kf(y)J(z,dy)\quad\quad\mbox{and}\quad\quad u(t,x)=1_\Omega(x)\left(f1_{(t_3,t_4]}\right)_{t_1,\Omega}(t,x).$$
	Note that for all $z\in\Omega$,
	$$u(t,x)=\begin{cases}0,&\mbox{ for all}\quad t_1<t\le t_3;\\
		2\int_{t_3}^{t\wedge t_4}\int_\Omega\int_Kf(y)J(z,dy)p_{t-s}^\Omega(x,dz)ds=2\int_{(t-t_4)_+}^{t-t_3}P_s^\Omega T_f(x)ds,&\mbox{ for all}\quad t>t_3.\end{cases}$$
	Therefore, for all $t>t_3$ and $\sigma>0$,
	\begin{align*}
		\mathcal{E}_\sigma(u,u)&=4\int_{(t-t_4)_+}^{t-t_3}\int_{(t-t_4)_+}^{t-t_3}\mathcal{E}_\sigma\left(P_s^\Omega T_f,P_{s'}^\Omega T_f\right)ds'ds=\frac{4}{\sigma}\int_{(t-t_4)_+}^{t-t_3}\int_{(t-t_4)_+}^{t-t_3}\left(P_s^\Omega T_f-P_{s+\sigma}^\Omega T_f,P_{s'}^\Omega T_f\right)ds'ds\\
		&=-\frac{4}{\sigma}\int_{(t-t_4)_+}^{t-t_3}\int_{(t-t_4)_+}^{t-t_3}\int_0^\sigma\frac{d}{d\sigma'}\left(P_{s+\sigma'}^\Omega T_f,P_{s'}^\Omega T_f\right)d\sigma'ds'ds\\
		&=-\frac{4}{\sigma}\int_{(t-t_4)_+}^{t-t_3}\int_{(t-t_4)_+}^{t-t_3}\int_0^\sigma\frac{d}{ds}\left(P_{s+\sigma'}^\Omega T_f,P_{s'}^\Omega T_f\right)d\sigma'ds'ds\\
		&=\frac{4}{\sigma}\int_{(t-t_4)_+}^{t-t_3}\int_0^\sigma\left(P_{(t-t_4)_++\sigma'}^\Omega T_f-P_{t-t_3+\sigma'}^\Omega T_f,P_{s'}^\Omega T_f\right)d\sigma'ds'\le\frac{4}{\sigma}\int_{(t-t_4)_+}^{t-t_3}\int_0^\sigma\left(P_{(t-t_4)_+}^\Omega T_f,P_{s'+\sigma'}^\Omega T_f\right)d\sigma'ds'\\
		&=\frac{4}{\sigma}\int_0^\sigma\int_\Omega\int_\Omega\int_K\left(\int_{(t-t_4)_+}^{t-t_3}\int_\Omega\int_Kf(z)J(y,dz)p_{s'+\sigma'}^\Omega(x,dy)ds'\right)f(z')J(x',dz')p_{(t-t_4)_+}^\Omega(x,dx')d\mu(x)d\sigma'\\
		&\le\frac{2\|f\|_\infty^2}{\sigma}\int_0^\sigma\int_\Omega\int_\Omega\int_KJ(x',dz')p_{(t-t_4)_+}^\Omega(x',dx)d\mu(x')d\sigma'\\
		&\le 2\|f\|_\infty^2\int_\Omega\int_KJ(x',dz')d\mu(x')=2\|f\|_\infty^2\int_\Omega\int_K(\phi(x')-\phi(z'))^2J(x',dz')d\mu(x')\le 2\|f\|_\infty^2\mathcal{E}(\phi)<\infty.
	\end{align*}
	Here we use (\ref{levy+}) in the sixth line and $\int_\Omega p_t^\Omega(x',dx)\le 1$ in the seventh line.
	
	Same as in (iii) of Lemma \ref{ext-F}, we see $u(t,\cdot)\in\mathcal{F}'(\Omega)$ and is caloric on $(t_1,\infty)\times\Omega$.
\end{proof}

\subsection{The tail in $(\mathrm{PHI})$}

\begin{proposition}\label{phi+t}
	Assume that $(\mathrm{PHI}^0)$ holds. Then, there exist $\delta_\ast>0$ and $\lambda_\ast\in(0,1)$ such that given any $D_\ast\ge\delta_\ast$, there exists $C_\ast=C_\ast(D_\ast)>0$ such that for every ball $B=B(x_0,R)$ with $0<R<\overline{R}$, all $t_4<t_3<t_2<t_1<t_0$ satisfying
	\begin{equation}\label{intrv}
		\delta_\ast W(B)\le t_i-t_{i+1}\le D_\ast W(B)\quad\quad\mbox{for all}\quad\quad i=0,1,2,3
	\end{equation}
	and all $u:(t_4,t_0]\to\mathcal{F}'(B)$ that is non-negative, bounded and caloric on $Q:=(t_4,t_0]\times B$,
	\begin{equation}\label{phi'var}
		\mathop{\mathrm{esup}}_{(t_3,t_2]\times\lambda_\ast B}u\le C_\ast\left(\mathop{\mathrm{einf}}_{(t_1,t_0]\times\lambda_\ast B}u+T_Q^{(3/4;1)}(u_-)\right).
	\end{equation}
\end{proposition}

In particular, $(\mathrm{PHI}^+)$, $(\mathrm{PHI}^0)$ and all presentations of the parabolic Harnack inequality investigated by \cite{ckw-p} are equivalent with each other.

\begin{proof}
	It suffices to consider the case $T_Q^{(3/4;1)}(u_-)<\infty$. Let
	$$\delta_\ast=C_W\left(\delta_2\vee(\delta_4-\delta_2)\right)(\lambda_1^{\beta_2}\vee\lambda_1^{\beta_1}).$$
	For any $D_\ast\ge\delta_\ast$, any ball $B=B(x_0,R)$, any $t_0,\dotsc,t_4$ satisfying (\ref{intrv}) and desired $u$, define
	\begin{equation}\label{uep}
		u_\varepsilon(t,x):=\frac{1}{\varepsilon}\int_0^\varepsilon u(t-s',x)ds'\quad\quad\mbox{and}\quad\quad \hat{u}_\varepsilon:=u_\varepsilon+((u_\varepsilon)_-)_{t_4+\varepsilon,\frac{3}{4}B}
	\end{equation}
	for any $0<\varepsilon<\delta_\ast W(B)-\delta_2W(\lambda_1B)$. Clearly $u_\varepsilon$ is caloric and bounded on $Q_\varepsilon=(t_4+\varepsilon,t_0]\times B$. Further,
	$$T_{(t_4+\varepsilon,t_0]\times B}^{(3/4)}((u_\varepsilon)_-)\le\sup\limits_{t_4+\varepsilon<t\le t_0}\frac{1}{\varepsilon}\int_{t-\varepsilon}^t T_B^{(3/4)}(u_-(s,\cdot))ds\le\frac{1}{\varepsilon}T_Q^{(3/4;1)}(u_-)<\infty.$$
	Thus by Lemma \ref{ext-F}, $\hat{u}_\varepsilon$ is caloric and bounded on $(t_4+\varepsilon,t_0]\times\frac{3}{4}B$. It is easy to check that $\hat{u}_\varepsilon\ge 0$ globally with the help of Proposition \ref{gpmp}.
	
	By the definition of esup, einf in (\ref{e'sup}), there exist $s_{u,\varepsilon}\in(t_3,t_2]$ and $s'_{u,\varepsilon}\in(t_1,t_0]$ such that
	$$\mathop{\mathrm{esup}}\limits_{(t_3,t_2]\times\lambda_0B}u_\varepsilon=\mathop{\mathrm{esup}}\limits_{\lambda_0B}u_\varepsilon(s_{u,\varepsilon},\cdot),\quad\quad\mathop{\mathrm{einf}}\limits_{(t_1,t_0]\times\lambda_0B}u_\varepsilon=\mathop{\mathrm{einf}}\limits_{\lambda_0B}u_\varepsilon(s'_{u,\varepsilon},\cdot).$$
	Define
	$$n_{u,\varepsilon}:=\left\lceil\frac{s'_{u,\varepsilon}-s_{u,\varepsilon}}{(\delta_4-\delta_2)W(\lambda_1B)}\right\rceil\quad\quad\mbox{and}\quad\quad r_{u,\varepsilon}:=\frac{3}{4\lambda_1}W^{-1}\left(x_0,\frac{s'_{u,\varepsilon}-s_{u,\varepsilon}}{n_{u,\varepsilon}(\delta_4-\delta_2)}\right).$$
	It is easy to prove that $n_{u,\varepsilon}>1$ (and hence $\frac{3}{4\lambda_1}W^{-1}(x_0,\frac{1}{2}W(\lambda_1B))\le r_{u,\varepsilon}\le\frac{3}{4}R$), while on the other hand,
	$$n_{u,\varepsilon}\le 1+\frac{3D_\ast W(B)}{(\delta_4-\delta_2)W(\lambda_1B)}\le 1+\frac{3C_WD_\ast}{\delta_4-\delta_2}(\lambda_1^{-\beta_2}\vee\lambda_1^{-\beta_1})=:N_\delta.$$
	Take $B_{u,\varepsilon}=B(x_0,r_{u,\varepsilon})$ and
	$$\tau_n=s_{u,\varepsilon}+(n\delta_4-(n+1)\delta_2)W(x_0,\lambda_1r_{u,\varepsilon})\quad\quad\mbox{for}\quad\quad n=0,\dotsc,n_{u,\varepsilon}-1.$$
	Note that $\tau_0>\varepsilon+t_4$ by the selection of $\delta_\ast$ and $\varepsilon$. Take further
	$$\lambda_\ast=\frac{3}{4}(2C_W)^{-1/\beta_1}\lambda_0\quad\quad\mbox{so that}\quad\quad\lambda_\ast R\le\lambda_0r_{u,\varepsilon}.$$
	Applying $(\mathrm{PHI}^0)$ to $(t,x)\mapsto\hat{u}_\varepsilon(t+\tau_n,x)$ on $(0,\delta_4W(\lambda_1B_{u,\varepsilon})]\times B_{u,\varepsilon}$ successively, we obtain
	\begin{align*}
		\mathop{\mathrm{esup}}\limits_{(t_3,t_2]\times\lambda_\ast B}u_\varepsilon\le&\ \mathop{\mathrm{esup}}\limits_{(t_3,t_2]\times\lambda_0B_{u,\varepsilon}}u_\varepsilon\le\mathop{\mathrm{esup}}\limits_{(t_3,t_2]\times\lambda_0B}u_\varepsilon=\mathop{\mathrm{esup}}_{\lambda_0B_{u,\varepsilon}}u_\varepsilon(s_{u,\varepsilon},\cdot)\le\mathop{\mathrm{esup}}_{(t,x)\in(\delta_1W(\lambda_1B_{u,\varepsilon}),\delta_2W(\lambda_1B_{u,\varepsilon})]\times\lambda_0B_{u,\varepsilon}}\hat{u}_\varepsilon(\tau_0+t,x)\\
		\le&\ C\mathop{\mathrm{einf}}_{(t,x)\in(\delta_3W(\lambda_1B_{u,\varepsilon}),\delta_4W(\lambda_1B_{u,\varepsilon})]\times\lambda_0B_{u,\varepsilon}}\hat{u}_\varepsilon(\tau_0+t,x)\\
		\le&\ C\mathop{\mathrm{einf}}_{\lambda_0B_{u,\varepsilon}}\hat{u}_\varepsilon(\tau_0+\delta_4W(\lambda_1B_{u,\varepsilon}),\cdot)=C\mathop{\mathrm{einf}}_{\lambda_0B_{u,\varepsilon}}\hat{u}_\varepsilon(\tau_1+\delta_2W(\lambda_1B_{u,\varepsilon}),\cdot)\\
		\le&\ C\mathop{\mathrm{esup}}_{(t,x)\in(\delta_1W(\lambda_1B_{u,\varepsilon}),\delta_2W(\lambda_1B_{u,\varepsilon})]\times\lambda_0B_{u,\varepsilon}}\hat{u}_\varepsilon(\tau_1+t,x)\le C\mathop{\mathrm{einf}}_{\lambda_0B_{u,\varepsilon}}\hat{u}_\varepsilon(\tau_1+\delta_4W(\lambda_1B_{u,\varepsilon}),\cdot)\\
		\le&\ \cdots\le C^{n_u-1}\mathop{\mathrm{einf}}_{\lambda_0B_{u,\varepsilon}}\hat{u}_\varepsilon(\tau_{n_{u,\varepsilon}-1}+\delta_4W(\lambda_1B_{u,\varepsilon}),\cdot)=C^{n_{u,\varepsilon}-1}\mathop{\mathrm{einf}}_{\lambda_0B_{u,\varepsilon}}\hat{u}_\varepsilon(s'_{u,\varepsilon},\cdot)\\
		=&\ C^{n_{u,\varepsilon}-1}\mathop{\mathrm{einf}}\limits_{\lambda_0B_{u,\varepsilon}}\left(u_\varepsilon(s'_{u,\varepsilon},\cdot)+2\int_{t_4+\varepsilon}^{s'_{u,\varepsilon}}\int_{\frac{3}{4}B}\int_{B^c}(u_\varepsilon)_-(s,y)J(z,dy)p_{t-s}^{\frac{3}{4}B}(x,dz)ds\right)\\
		\le&\ C^{N_\delta-1}\mathop{\mathrm{einf}}\limits_{\lambda_0B_{u,\varepsilon}}u_\varepsilon(s'_{u,\varepsilon},\cdot)+2C^{N_\delta-1}\mathop{\mathrm{esup}}\limits_{\lambda_0B_{u,\varepsilon}}\int_{t_4+\varepsilon}^{s'_{u,\varepsilon}}\left\{\frac{1}{\varepsilon}\int_{s-\varepsilon}^sT_B^{(3/4)}(u_-(s',\cdot))ds'\int_{\frac{3}{4}B}p_{s'_{u,\varepsilon}-s}^{\frac{3}{4}B}(x,dz)\right\}ds\\
		\le&\ C^{N_\delta-1}\mathop{\mathrm{einf}}\limits_{(t_1,t_0]\times\lambda_\ast B}u_\varepsilon+2C^{N_\delta-1}T_Q^{(3/4;1)}(u_-).
	\end{align*}
	Letting $\varepsilon\downarrow 0$ and using the weak $L^2(B)$-continuity of $u$, we obtain
	$$\mathop{\mathrm{esup}}\limits_{(t_3,t_2]\times\lambda_\ast B}u\le 2C^{N_\delta-1}\left(\mathop{\mathrm{einf}}\limits_{(t_1,t_0]\times\lambda_\ast B}u+T_Q^{(3/4;1)}(u_-)\right),$$
	that is, (\ref{phi'var}) holds with $C_\ast=2C^{N_\delta-1}$.
\end{proof}

\section{Properties involving the upper jumping smoothness}\label{t-ujs}

\subsection{Proving $(\mathrm{UJS})$ from $(\mathrm{PHI}^0)$}\label{p'ujs}

\begin{proposition}\label{ujs}
	$(\mathrm{PHI}^0)\Rightarrow(\mathrm{UJS})$.
\end{proposition}

The proof here is purely analytic, whose idea is based on the stochastic proof of \cite[Proposition 3.3]{ckw-p}.

\begin{proof}
	To begin with, by $(\mathrm{PHI}^0)$, we see for every ball $B$, each $\delta_3W(\lambda_1B)<t\le\delta_4W(\lambda_1B)$ and every globally non-negative function $u$ that is bounded and caloric on $(0,\delta_4W(\lambda_1B)]\times B$,
	\begin{equation}\label{pmv+}
		\mathop{\mathrm{esup}}\limits_{(\delta_1W(\lambda_1B),\delta_2W(\lambda_1B)]\times\lambda_0B}u\le C\mathop{\mathrm{einf}}\limits_{(\delta_3W(\lambda_1B),\delta_4W(\lambda_1B)]\times\lambda_0B}u\le C\mathop{\mathrm{einf}}\limits_{\lambda_0B}u(t,\cdot)\le C\fint_{\lambda_0B}u(t,z)d\mu(z).
	\end{equation}
	
	Given any two different points $x_0,y_0\in M$. For any $0<r_0\le\frac{1}{6}d(x_0,y_0)$ and $0<r\le\frac{1}{2}d(x_0,y_0)$, the two balls $B=B(x_0,r)$ and $B_0=B(y_0,r_0)$ are disjoint with $d(B,B_0)\ge 2r_0$. For any $0<h<\delta_2W(\lambda_1B)$, let
	$$f_h(t,z)=\frac{1}{h}\cdot 1_{(\delta_2W(\lambda_1B)-h,\delta_2W(\lambda_1B)]}(t)1_{B_0}(z).$$
	By Proposition \ref{cal+}, $u_h:=(f_h)_{0,B}$ is strongly caloric on $(0,\delta_4W(\lambda_1B)]\times B$. Clearly $u_h$ is non-negative on $M$. Thus by applying (\ref{pmv+}) with $t=\delta_4W(\lambda_1B)$ and the absolute continuity of the semigroup $P_t^B$, for $\mu$-a.e.\ $x\in\lambda_0B$,
	\begin{align*}
		\int_{B_0}J(x,dy)=&\ \lim\limits_{h\to 0}u_h(\delta_2W(\lambda_1B),x)\le\limsup\limits_{h\to 0}\frac{C}{\mu\left(\lambda_0B\right)}\int_{\lambda_0B}u_h(\delta_4W(\lambda_1B),z)dz\\
		=&\ \frac{C}{\mu\left(\lambda_0B\right)}\int_{\lambda_0B}P_{\delta_2W(\lambda_1B)}^BJ(\cdot,B_0)(z)dz\le\frac{C}{\mu\left(\lambda_0B\right)}\int_{\lambda_0B}\int_{B_0}J(z,dy)d\mu(z).
	\end{align*}
	Therefore, for all $\varepsilon<\lambda_0r$ and $\mu$-a.e.\ $x\in B\left(x_0,\lambda_0r-\varepsilon\right)$,
	$$\int_{B(y_0,r_0)}\fint_{B(x,\varepsilon)}J(y,dx)d\mu(y)=\fint_{B(x,\varepsilon)}\int_{B(y_0,r_0)}J(x,dy)d\mu(x)\le\frac{C}{\mu\left(\lambda_0B\right)}\int_{\lambda_0B}\int_{B(y_0,r_0)}J(z,dy)d\mu(z).$$
	Since $y_0,r_0$ are arbitrary, we have for $\mu$-a.e.\ $y\in B^c$,
	$$\fint_{B(x,\varepsilon)}J(y,dx)\le\frac{C}{\mu\left(\lambda_0B\right)}\int_{\lambda_0B}J(y,dz).$$
	Thus $J(y,dx)$ is absolutely continuous and the density $x\mapsto J(y,x)$ satisfies
	$$J(y,x)=\lim\limits_{\varepsilon\to 0}\fint_{B(x,\varepsilon)}J(y,dx)\le\frac{C}{\mu\left(\lambda_0B\right)}\int_{\lambda_0B}J(y,dz)=\frac{C}{\mu\left(\lambda_0B\right)}\int_{\lambda_0B}J(y,z)d\mu(z)$$
	for $\mu$-a.e.\ $x\in \lambda_0B$. $(\mathrm{UJS})$ follows by the standard covering technique.
\end{proof}

\subsection{Regularity improved by $(\mathrm{UJS})$}\label{uj1}

The following lemma is simple but fundamental:

\begin{lemma}\label{L2-J}
	$(\mathrm{cap}_\le)+(\mathrm{UJS})\Rightarrow(\mathrm{TJ})$, and hence $(\mathrm{J}_\le)$ also holds.
\end{lemma}

\begin{proof}
	Fix an arbitrary $B=B(x_0,R)$ and $0<\lambda<1$. Let $\lambda'=\frac{1+\lambda}{2}$. By $(\mathrm{cap}_\le)$, there exists $\phi\in\mathrm{cutoff}(\lambda'B,B)$ such that $\mathcal{E}(\phi,\phi)\le C_{\lambda'}\frac{\mu(B)}{W(B)}$. Therefore,
	$$j\left(\lambda'B\times B^c\right)=\int_{\lambda'B\times B^c}dj(x,y)=\int_{\lambda'B\times B^c}(\phi(x)-\phi(y))^2dj(x,y)\le\mathcal{E}(\phi,\phi)\le C_{\lambda'}\frac{\mu(B)}{W(B)}.$$
	For all $x\in\lambda B$, $z\in B\left(x,\frac{1-\lambda}{2}R\right)$ (so that $z\in\lambda'B$) and $y\in B^c$, since $d(x,z)<\frac{1-\lambda}{2}R\le\frac{1}{2}d(x,y)$, then by $(\mathrm{UJS})$,
	$$\int_{B^c}J(x,dy)=\int_{B^c}J(x,y)d\mu(y)\le\int_{B^c}C_{\mathrm{UJS}}\fint_{B(x,\frac{1-\lambda}{2}R)}J(z,y)d\mu(z)d\mu(y)\le\frac{C(\lambda)}{\mu(B)}\int_{\lambda'B\times B^c}dj(z,y)\le\frac{C(\lambda)C_{\lambda'}}{W(B)},$$
	which is exactly $(\mathrm{TJ})$. A similar treatment on $y$ yields $(\mathrm{TJ})+(\mathrm{UJS})\Rightarrow(\mathrm{J}_\le)$.
\end{proof}

\begin{lemma}\label{T+byT-}
	If $(\mathrm{UJS})$ holds, then there exists $C'>0$ such that for every $u$ that is non-negative and caloric on $Q=(t_1,t_2]\times B$, for each $\phi\in\mathrm{cutoff}\left(\frac{17}{24}B,\frac{3}{4}B\right)$,
	\begin{equation}\label{up-t+}
		T_Q^{(2/3;1)}(u_+)\le C'\left\{\sup\limits_{t_1<s\le t_2}\|u(s,\cdot)\|_{L^\infty(B)}\left(1+\frac{t_2-t_1}{\mu(B)}\mathcal{E}(\phi,\phi)\right)+T_Q^{(3/4;1)}(u_-)\right\}.
	\end{equation}
\end{lemma}

In particular, if $(\mathrm{cap}_\le)$ holds, then
\begin{equation}\label{up-t-}
	T_Q^{(2/3;1)}(u_+)\le C'\left\{\sup\limits_{t_1<s\le t_2}\|u(s,\cdot)\|_{L^\infty(B)}\left(1+\frac{t_2-t_1}{W(B)}\right)+T_Q^{(3/4;1)}(u_-)\right\}.
\end{equation}

\begin{proof}
	Clearly (\ref{up-t+}) holds automatically if $\|u\|_Q=\infty$. For $0\le\|u\|_Q<\infty$, with an arbitrary $\varepsilon>0$ we take
	$$u'(t,x):=(\|u\|_Q\vee\varepsilon)^{-1}u(t,x).$$
	Then $0\le u'\le 1$ on $Q$. Since $u'(t,\cdot)\in\mathcal{F}'(B)$, by decomposing $u'$ into $w+(u'-w)$, where $w\in\mathcal{F}$ and $\mathrm{supp}(u'-w)\subset B^c$, we see by \cite[Lemma 3.7]{ghh18} that for every $0<t\le T$ and $\phi\in\mathrm{cutoff}\left(\frac{17}{24}B,\frac{3}{4}B\right)$,
	$$\mathcal{E}^{(J)}\left(u',\frac{\phi^2}{u'+2}\right)\le 3\mathcal{E}^{(J)}(\phi,\phi)-2\int_{B\times B^c}\phi^2(x)\frac{u'(t,y)+2}{u'(t,x)+2}dj(x,y),$$
	where $\mathcal{E}^{(J)}(u,v)$ stands for the second term in (\ref{bd}). Meanwhile, by the chain rule and Leibniz's rule of $\Gamma^{(L)}$,
	\begin{align*}
		\mathcal{E}^{(L)}\left(u',\frac{\phi^2}{u'+2}\right)=&\int_Bd\Gamma^{(L)}\left(u',\frac{\phi^2}{u'+2}\right)=\int_B\frac{2\phi}{u'+2}d\Gamma^{(L)}(u',\phi)-\frac{\phi^2}{(u'+2)^2}d\Gamma^{(L)}(u',u')\\
		\le&\int_B\frac{\phi^2}{(u'+2)^2}d\Gamma^{(L)}(u',u')+\int_Bd\Gamma^{(L)}(\phi,\phi)-\int_B\frac{\phi^2}{(u'+2)^2}d\Gamma^{(L)}(u',u')=\mathcal{E}^{(L)}(\phi,\phi).
	\end{align*}
	It follows that
	\begin{align}
		-\frac{d}{dt}\int_B\phi^2&\log(u'+2)d\mu=\mathcal{E}\left(u',\frac{\phi^2}{u'+2}\right)=\mathcal{E}^{(J)}\left(u',\frac{\phi^2}{u'+2}\right)+\mathcal{E}^{(L)}\left(u',\frac{\phi^2}{u'+2}\right)\notag\\
		\le&\ 3\mathcal{E}^{(J)}(\phi,\phi)-2\int_{B\times B^c}\phi^2(x)\frac{u'(t,y)+2}{u'(t,x)+2}dj(x,y)+\mathcal{E}^{(L)}(\phi,\phi)\notag\\
		\le&\ 3\mathcal{E}(\phi,\phi)-2\int_{B\times\{y:u(t,y)<0\}}\phi^2(x)\frac{u'(t,y)+2}{u'(t,x)+2}dj(x,y)-2\int_{B\times(B^c\cap\{y:u(t,y)\ge 0\})}\phi^2(x)\frac{u'(t,y)+2}{u'(t,x)+2}dj(x,y)\notag\\
		\le&\ 3\mathcal{E}(\phi,\phi)+\int_{\frac{3}{4}B\times B^c}u'_-(t,y)dj(x,y)-2\int_{\frac{17}{24}B\times B^c}\frac{u'_+(t,y)}{3}dj(x,y)\notag\\
		\le&\ 3\mathcal{E}(\phi,\phi)+\frac{1}{\|u\|_Q\vee\varepsilon}\mu\left(\frac{3}{4}B\right)T_B^{(3/4)}(u_-(t,\cdot))-\frac{2/3}{\|u\|_Q\vee\varepsilon}\int_{\frac{17}{24}B\times B^c}u_+(t,y)dj(x,y),\label{log-T}
	\end{align}
	where we used the following facts for the fourth line:
	\begin{align*}
		2\le u'(t,x)+2\le 3,&\quad\quad\mbox{for all}\quad(t,x)\in Q;\\
		-(u'(t,y)+2)=u'_-(t,y)-2<u'_-(t,y),&\quad\quad\mbox{whenever}\quad u(t,y)<0;\\
		u'(t,y)+2>u'(t,y)=u'_+(t,y),&\quad\quad\mbox{whenever}\quad u(t,y)\ge 0.
	\end{align*}
	Note that for all $x\in B\left(x_0,\frac{2}{3}R\right)$, $B\left(x,\frac{1}{24}R\right)\subset B\left(x_0,\frac{17}{24}R\right)$ and $B(x_0,R)\subset B\left(x,\frac{5}{3}R\right)$. Consequently, by $(\mathrm{UJS})$,
	\begin{align*}
		T_Q^{(\frac{2}{3};1)}(u_+)=\ &\int_{t_1}^{t_2}\sup\limits_{x\in\frac{2}{3}B}\int_{B^c}u_+(t,y)J(x,y)d\mu(y)\le C\int_{t_1}^{t_2}\sup\limits_{x\in\frac{2}{3}B}\int_{B^c}\fint_{B\left(x,\frac{R}{24}\right)}u_+(t,y)J(z,y)d\mu(z)d\mu(y)\\
		\le\ &\frac{C^2}{\mu(B)}\int_{t_1}^{t_2}\int_{B^c}\int_{\frac{17}{24}B}u_+(t,y)J(z,y)d\mu(z)d\mu(y)\\
		\le\ &\frac{C^2(\|u\|_Q\vee\varepsilon)}{\mu(B)}\int_{t_1}^{t_2}\left\{\frac{d}{dt}\int_B\phi^2\log(u'(t,\cdot)+2)d\mu+3\mathcal{E}(\phi,\phi)+\frac{\mu\left(\frac{3}{4}B\right)}{\|u\|_Q\vee\varepsilon}T_B^{(3/4)}(u_-(t,\cdot))\right\}dt\\
		\le\ &\frac{C^2(\|u\|_Q\vee\varepsilon)}{\mu(B)}\left\{\left[\int_B\phi^2\log(u'(t,\cdot)+2)d\mu\right]\bigg|_{t=t_1}^{t_2}+3(t_2-t_1)\mathcal{E}(\phi,\phi)\right\}+C^2\int_{t_1}^{t_2}T_B^{(3/4)}(u_-(t,\cdot))dt\\
		\le\ &C^2(\|u\|_Q\vee\varepsilon)\left\{\log 3+\frac{3(t_2-t_1)}{\mu(B)}\mathcal{E}(\phi,\phi)\right\}+C^2T_Q^{(3/4;1)}(u_-),
	\end{align*}
	where we used (\ref{log-T}) for the third line. Thus (\ref{up-t+}) follows with $C'=3C^2$ by letting $\varepsilon\to 0$.
\end{proof}

As a consequence, we find:

\begin{proposition}\label{ph+}
	$(\mathrm{PHI})\Leftrightarrow(\mathrm{PHI}^+)$.
\end{proposition}

\begin{proof}
	Trivially ``$\Leftarrow$'' holds. On the converse, assume that $(\mathrm{PHI})$ holds. Then, by Proposition \ref{ujs} we see $(\mathrm{UJS})$ holds; while by \cite[Proposition 3.9]{L0} and the known implication $(\mathrm{LLE})\Rightarrow(\mathrm{cap}_\le)$ we see $(\mathrm{cap}_\le)$ holds. Hence by (\ref{up-t-}),
	$$T_{Q_-}^{(2/3;1)}(u_+)\le C\left(\sup\limits_{Q_-}u+T_{Q_-}^{(3/4;1)}(u_-)\right)\le C\left(\sup\limits_{Q_-}u+T_Q^{(3/4;1)}(u_-)\right).$$
	Combining $(\mathrm{PHI})$, we obtain
	$$\sup\limits_{Q_-}u+T_{Q_-}^{(2/3;1)}(u_+)\le(1+C'')\sup\limits_{Q_-}u+C''T_Q^{(3/4;1)}(u_-)\le K(1+C'')\inf\limits_{Q_+}u+\left(C''+2\delta K(1+C'')\right)T_Q^{(3/4;1)}(u_-),$$
	showing the inequality (\ref{phi+}). The proof is then completed.
\end{proof}

Based on these results, we obtain the following important result:

\begin{proposition}\label{sp-t}
	Assume that $(\mathrm{Gcap})$, $(\mathrm{PI})$ and $(\mathrm{UJS})$ hold, then given any cylinder $Q=(t_1,t_2]\times B$ (where $B$ is a metric ball), for every function $u:(t_1,t_2]\times\mathcal{F}'(B)$ that is non-negative, bounded and caloric on $Q$, if $T_Q^{(\lambda_1;1)}(u_-)<\infty$ for some $\lambda_1\in(0,1)$, then $u$ is uniformly continuous on $(t'_1,t_2]\times\lambda B$ for every $\lambda\in(0,\lambda_1)$.
\end{proposition}

\begin{proof}
	By Lemma \ref{T+byT-}, \cite[Proposition 6.2]{L0} and covering techniques.
\end{proof}

\subsection{Tail estimates improved by $(\mathrm{UJS})$}

\begin{lemma}\label{4-1}
	Assume that $(\mathrm{LLE})$ and $(\mathrm{UJS})$ hold. Let $\delta,\delta'\in(0,\frac{1}{4})$ be any two constants such that
	\begin{equation}\label{del-0}
		4\delta W(x,r)\le W(x,\delta_Lr)\quad\quad\mbox{and}\quad\quad W(x,2\delta'r)\le\delta W(x,\delta_Lr)
	\end{equation}
	for all $x\in M$ and $r>0$, where $\delta_L$ comes from $(\mathrm{LLE})$. Then, there exists $C_1=C_1(\delta,\delta')>0$ such that for every ball $B=B(x_0,r)$, all $t_0\in\mathbb{R}$ and non-negative $\hat{\mu}$-measurable function $f:(t_0-4\delta W(B),t_0]\to L_{loc}^1(M)$ with $\mathrm{supp}(f(t,\cdot))\subset(2B)^c$ for all $t_0-4\delta W(B)<t\le t_0$,
	$$\mathop{\mathrm{esup}}\limits_{Q_-}f_{t_0-4\delta W(B),2\delta'B}\le C_1\mathop{\mathrm{einf}}\limits_{Q_+}f_{t_0-4\delta W(B),B},$$
	where $Q_-:=(t_0-3\delta W(B),t_0-2\delta W(B)]\times\delta'B$ and $Q_+:=(t_0-\delta W(B),t_0]\times\delta'B$.
\end{lemma}

The proof is essentially contained in \cite[Lemma 5.3]{ckk}, but we reduce the condition $(\mathrm{UE})$.

\begin{proof}
	Denote $t_k=t_0-k\delta W(B)$ for $i=1,\dotsc,4$. Note that for all $z\in\frac{3}{2}\delta'B$, $z'\in B(z,\frac{1}{2}\delta'r)$ and $y\in(2B)^c$,
	$$d(x_0,z')\le d(x_0,z)+d(z,z')<\frac{3}{2}\delta'r+\frac{1}{2}\delta'r=2\delta'r,\quad\quad d(z,z')<\frac{1}{2}\delta'r\le\frac{1}{2}\left(2r-\frac{3}{2}\delta'r\right)\le\frac{1}{2}d(z,y).$$
	Therefore, by $(\mathrm{UJS})$ and $(\mathrm{UE})$, for all $t_3<t\le t_2$ and $\mu$-a.e. $x\in\delta'B$ (that is, $(t,x)\in Q_-$),
	\begin{align*}
		f_{t_4,2\delta'B}(t,x)=&\ 2\int_{t_4}^t\int_{(2B)^c}f(s,y)\int_{2\delta'B}p_{t-s}^{2\delta'B}(x,z)J(y,z)d\mu(z)d\mu(y)ds\\
		\le&\ 2\int_{t_4}^t\int_{(2B)^c}f(s,y)\int_{2\delta'B}p_{t-s}^{2\delta'B}(x,z)\left(\frac{C_{\mathrm{UJS}}}{V\left(z,\frac{1}{2}\delta'r\right)}\int_{B\left(z,\frac{1}{2}\delta'r\right)}J(z',y)d\mu(z')\right)d\mu(z)d\mu(y)ds\\
		\le&\ 2\int_{t_4}^t\int_{(2B)^c}f(s,y)\frac{C_{\mathrm{UJS}}C_{\mathrm{VD}}3^\alpha}{\mu(\delta'B)}\int_{3\delta'B}J(z',y)d\mu(z')d\mu(y)ds\\
		\le&\ \frac{C(\delta')}{\mu(B)}\int_{t_4}^{t_2}\int_{(2B)^c}f(s,y)\int_{3\delta'B}J(y,dz)d\mu(y)ds.
	\end{align*}
	On the other hand, given any $t_1<t'\le t_0$, $\mu$-a.e. $x'_1\in\delta'B$ (that is, $(t',x'_1)\in Q_+$) and $x'_2\in 2\delta'B$, note that
	\begin{align*}
		t'-s<&\ t_0-t_4=4\delta(B)\le W(\delta_LB);\\
		d(x_0,x'_i)<&\ 2\delta'r=\delta_LW^{-1}(x_0,W(2\delta'\delta_L^{-1}B))\le\delta_LW^{-1}(x_0,\delta W(B))\le\delta_LW^{-1}(x_0,t_1-t_2)\le\delta_LW^{-1}(x_0,t'-s)
	\end{align*}
	for all $t_4<s\le t_2$. Thus by $(\mathrm{LLE})$,
	\begin{align*}
		f_{t_4,B}(t',x'_1)\ge&\ 2\int_{t_4}^{t_2}\int_{(2B)^c}f(s,y)\int_{3\delta'B}p_{t'-s}^B(x'_1,x'_2)J(y,dx'_2)d\mu(y)ds\\
		\ge&\ 2\int_{t_4}^{t_2}\int_{(2B)^c}f(s,y)\frac{\int_{3\delta'B}c_LJ(y,dx'_2)d\mu(y)}{V(x_0,W^{-1}(x_0,t'-s))}ds\ge 2\int_{t_0}^{t_2}\int_{(2B)^c}f(s,y)\frac{c_L\int_{3\delta'B}J(y,dx'_2)d\mu(y)}{V\left(x_0,W^{-1}(x_0,4\delta W(B))\right)}ds\\
		\ge&\ \frac{c(\delta)}{\mu(B)}\int_{t_4}^{t_2}\int_{(2B)^c}f(s,y)\int_{3\delta'B}J(y,dx'_2)d\mu(y)ds.
	\end{align*}
	Thus the proof is completed by taking $C_1=C(\delta')/c(\delta)$.
\end{proof}

\begin{proposition}\label{tail}
	Under the same conditions as in Lemma \ref{4-1}, with same $\delta$ and $\delta'$, there exists $\delta''\in(0,1)$ such that for any globally non-negative $u$ that is caloric on $(t_0-\delta W(B),t_0+3\delta W(\frac{1}{2}B)]\times B$,
	$$T_{(t_0-\delta W(\delta'B),t_0]\times B}^{(\delta'\delta'';1)}(u_+)\le C\mathop{\mathrm{einf}}\limits_{(t_0+2\delta W(\frac{1}{2}B),t_0+3\delta W(\frac{1}{2}B)]\times\frac{1}{2}\delta'B}u.$$
\end{proposition}

\begin{proof}
	Observe for all $x'\in\delta''\delta'B$, $z\in B^c$ and $x\in B(x',\delta''\delta'R)$ that
	$$d(x,x')<\delta''\delta'R<\frac{1}{4}R\le\frac{1}{2}d(z,x')\quad\mbox{and}\quad d(x,x_0)\le d(x,x')+d(x',x_0)<\delta''\delta'R+\delta''\delta'R=2\delta''\delta'R.$$
	Thus by $(\mathrm{UJS})$ and $(\mathrm{VD})$,
	$$J(x',z)\le\frac{C}{V(x',\delta''\delta'R)}\int_{B(x',\delta''\delta'R)}J(x,z)d\mu(x)\le\frac{C^2}{\mu(\delta''\delta'B)}\int_{2\delta''\delta'B}J(x,z)d\mu(x).$$
	Further, take
	$$t'=t_0+\delta W(\delta'B),\quad\quad t''=t_0+3\delta W(\frac{1}{2}B)\quad\quad\mbox{and}\quad\quad\delta''\le\frac{\delta_L}{2}\left(\frac{\delta}{C_W}\right)^{1/\beta_1}.$$
	By the definition of $\delta$
	, for all $t_0-\delta W(\delta'B)<s\le t_0$ and $x\in 2\delta''\delta'B$,
	\begin{align*}
		t'-s&<(t_0+\delta W(\delta'B))-(t_0-\delta W(\delta'B))=2\delta W(x_0,\delta'R)\le W(x_0,\delta_L\delta'R);\\
		d(x_0,x)&<2\delta''\delta'R=2\delta''W^{-1}(x_0,W(\delta'B))\le\delta_LW^{-1}(x_0,\delta W(\delta'B))\le\delta_LW^{-1}(x_0,t'-s).
	\end{align*}
	Therefore, by $(\mathrm{LLE})$,
	$$p_{t'-s}^{\delta'B}(x_0,x)\ge\frac{c_L}{V\left(x_0,W^{-1}(x_0,t'-s)\right)}\ge\frac{c_L}{V\left(x_0,W^{-1}(x_0,2\delta W(\delta'B))\right)}\ge\frac{1}{D_1(\delta)\mu(\delta'B)}.$$
	Combining the inequalities above, we obtain
	\begin{align*}
		T_{(t_0-\delta W(\delta'B),t_0]\times B}^{(\delta'\delta'';1)}(u_+)=&\ 2\int_{t_0-\delta W(\delta'B)}^{t_0}\mathop{\mathrm{esup}}\limits_{x'\in\delta''\delta'B}\int_{B^c}u_+(s,z)J(x',z)d\mu(z)ds\\
		\le&\ 2\int_{t_0-\delta W(\delta'B)}^{t_0}\int_{B^c}\frac{C^2u_+(s,z)}{\mu(\delta''\delta'B)}\int_{2\delta''\delta'B}J(z,x)d\mu(x)d\mu(z)ds\notag\\
		\le&\ 2C^2\frac{D_1(\delta)\mu(\delta'B)}{\mu(\delta''\delta'B)}\int_{t_0-\delta W(\delta'B)}^{t_0}\int_{B^c}u_+(s,z)\int_{2\delta''\delta'B}p_{t'-s}^{\delta'B}(x_0,x)J(x,z)d\mu(x)d\mu(z)ds\notag\\
		\le&\ \frac{C'D_1(\delta)}{(\delta'')^\alpha}\cdot 2\int_{t_0-\delta W(\frac{1}{2}B)}^{t_0}\int_{B^c}1_{B^c}(z)u_+(s,z)\int_{\delta'B}p_{t'-s}^{\delta'B}(x_0,x)J(x,z)d\mu(x)d\mu(z)ds\notag\\
		=&\ \frac{C'D_1(\delta)}{(\delta'')^\alpha}\left(1_{B^c}u_+\right)_{t_0-\delta W(\frac{1}{2}B),\delta'B}(t',x_0)\le \frac{C'D_1(\delta)}{(\delta'')^\alpha}\mathop{\mathrm{esup}}\limits_{(t_0,t_0+\delta W(\frac{1}{2}B)]\times\frac{1}{2}\delta'B}\left(1_{B^c}u_+\right)_{t_0-\delta W(\frac{1}{2}B),\delta'B}\notag\\
		\le&\ D'_1(\delta,\delta'')\mathop{\mathrm{einf}}\limits_{(t_0+2\delta W(\frac{1}{2}B),t_0+3\delta W(\frac{1}{2}B)]\times\frac{1}{2}\delta'B}\left(1_{B^c}u_+\right)_{t_0-\delta W(\frac{1}{2}B),\frac{1}{2}B}\\
		\le&\ D'_1(\delta,\delta'')\mathop{\mathrm{einf}}\limits_{(t_0+2\delta W(\frac{1}{2}B),t_0+3\delta W(\frac{1}{2}B)]\times\frac{1}{2}\delta'B}u,
	\end{align*}
	Here we used Lemma \ref{4-1} for the sixth line, and then used the maximum principle for the last line.
\end{proof}

\section{Proof of $(\mathrm{PHI}^0)$}\label{Str}

Here we provide two proofs to the strong parabolic Harnack inequality, following the ideas by Nash \cite{nas} and Moser \cite{mos-p} respectively.

\subsection{The Nash approach}

\begin{proposition}\label{main'}
	$(\mathrm{LLE})+(\mathrm{UJS})\Rightarrow(\mathrm{PHI}^0)$.
\end{proposition}

\begin{proof}
	(1) Let $C_W$ be the constant in (\ref{W}). Without loss of generality, we assume $c_L\le 1$. Let
	\begin{equation}\label{del-1}
		\delta=\frac{1}{4C_W}\left(\frac{\delta_L}{3}\right)^{\beta_2}\quad\quad\mbox{so that}\quad\quad 4\delta W(x,3r')\le W(y,\delta_Lr')\quad\mbox{for all}\quad d(x,y)<r'.
	\end{equation}
	With some $\lambda_0\in(0,\frac{1}{12})$ to be determined later, given an arbitrary ball $B=B(x_0,R_0)$, let
	$$Q_0:=(0,4\delta W(B)]\times B,\quad\quad Q_+:=(3\delta W(B),4\delta W(B)]\times\lambda_0B,\quad\quad Q_-:=(\delta W(B),2\delta W(B)]\times\lambda_0B.$$
	
	Given any $u:(0,4\delta W(B)]\to\mathcal{F}'(B)$ that is caloric, bounded on $Q_0$ and globally non-negative, then
	$$M[u]:=\mathop{\mathrm{einf}}\limits_{Q_+}u<\infty$$ 
	and $T_{Q_0}^{(3/4;1)}(|u|)<\infty$ by Lemma \ref{T+byT-}. Thus by Proposition \ref{sp-t}, $u$ is continuous on $Q:=(0,4\delta W(B)]\times\frac{1}{2}B$.
	
	Now we prove there exists $K$ (independent of $B$) such that $\sup_{Q_-}u\le KM[u]$ for all $u$ with $0<M[u]<\infty$, which can be extended to the case $M[u]=0$ in the following way: since $M[u]=0$, then $M[u+\varepsilon]=\varepsilon>0$. Thus
	$$\mathrm{sup}_{Q_-}u=-\varepsilon+\mathrm{sup}_{Q_-}(u+\varepsilon)\le-\varepsilon+KM[u+\varepsilon]=(K-1)\varepsilon.$$
	As $\varepsilon\downarrow 0$, we have $\sup_{Q_-}u\le 0$, that is, (\ref{phi+}) holds with any $C>0$.
	
	(2) Let us show there exist $K>2$ and $1<\Lambda<2$ (both independent of $B$ and $u$) such that if there exists $(t_1,y_1)\in Q_-$ with $u(t_1,y_1)\ge KM[u]$, then there also exists
	$$(t_{n+1},y_{n+1})\in\left(t_n-\delta W\left(y_n,\Lambda^{1-n}r\right),t_n\right]\times B\left(y_n,3\Lambda^{1-n}r\right)$$
	inductively for all $n\ge 1$ such that $u(t_{n+1},y_{n+1})\ge K\Lambda^{\alpha n}M[u]$, where
	\begin{equation}\label{r-1}
		r:=\lambda_HR_0\quad\quad\mbox{with}\quad\quad\lambda_H=\frac{1-\Lambda^{-1}\lambda_0}{3}\wedge\left(\frac{1-\Lambda^{-\beta_1}}{2C_W}\right)^{\frac{1}{\beta_1}}.
	\end{equation}
	
	Suppose we have found $(t_i,y_i)$ for $i=1,\dotsc,n$. If $n\ge 2$, then
	$$d(y_i,y_n)\le\sum\limits_{j=i+1}^nd(y_{j-1},y_j)<\sum\limits_{j=i+1}^n3\Lambda^{2-j}r\le\frac{3\Lambda^{1-i}}{1-\Lambda^{-1}}r,$$
	which implies that
	$$d(x_0,y_n)\le d(x_0,y_1)+d(y_1,y_n)<\lambda_0R_0+\frac{3r}{1-\Lambda^{-1}}\le 2\lambda_0R_0<\frac{1}{6}R_0,$$
	and
	$$\frac{t_1-t_n}{\delta W(B)}=\sum\limits_{i=1}^{n-1}\frac{t_i-t_{i+1}}{\delta W(x_0,R_0)}<\sum\limits_{i=1}^{n-1}\frac{W(y_i,\Lambda^{1-i}r)}{W(x_0,R_0)}\le C_W\sum\limits_{i=1}^{n-1}\left(\frac{\Lambda^{1-i}r}{R_0}\right)^{\beta_1}\le\frac{C_W}{1-\Lambda^{-\beta_1}}\left(\frac{r}{R_0}\right)^{\beta_1}<\frac{1}{2}.$$
	The same is trivially true if $n=1$. Hence
	$$(t_n,y_n)\in\left(\frac{1}{2}\delta W(B),2\delta W(B)\right]\times 2\lambda_0B\quad\quad\mbox{and}\quad\quad \lambda_0B\Subset B\left(y_n,\frac{1}{3}R_0\right)\subset\frac{1}{2}B.$$
	Thus existence of the sequence $(t_n,y_n)$ implies unboundedness of $u$ on $Q_0$, which is a contradiction completing the proof.
	
	(3) Take $\delta'=\frac{\delta_L}{2}\left(C_W^{-1}\delta\right)^{1/\beta_1}$ so that (\ref{del-0}) holds. For simplicity, we denote
	$$r'_n=\Lambda^{1-n}r,\quad\quad r_n=2\delta'r'_n,\quad\quad B'_n=B(y_n,r'_n)\quad\quad\mbox{and}\quad\quad B_n=B(y_n,r_n).$$
	For any $0<\delta''<\frac{1}{2}$ and $\phi_n\in\mathrm{cutoff}(2B'_n,3B'_n)$, define for all $t_n-\delta W(B'_n)<t\le 4\delta W(B)$ and $x\in M$ that
	$$u_n(t,x):=\sup\limits_{(t_n-\delta W(B'_n),t_n]\times 3B'_n}u-u(t,x)\phi_n(x)+A_n(t),$$
	where
	$$A_n(t):=2T_{(t_n-\delta W(B_n),t]\times 2B'_n}^{(\delta'\delta'';1)}(u_+)=2\int_{t_n-\delta W(B_n)}^t\mathop{\mathrm{esup}}\limits_{x'\in\delta''B_n}\int_{(2B'_n)^c}u_+(s,z)J(x',dz)ds.$$
	
	Clearly $u_n\ge 0$ globally; while by Lemma \ref{trunc} with $U=\delta''B_n$, $\Omega=2B'_n$ and $\phi=\phi_n$, we see $u_n$ is supercaloric on
	$$Q'_n:=(t_n-\delta W(B_n),4\delta W(B)]\times\delta''B_n.$$
	Define
	$$Q_n:=\left(t_n-3\delta_0W(\delta''B_n),t_n-2\delta_0W(\delta''B_n)\right]\times\delta''B_n,\quad\quad D_n:=\left\{(t,x)\in Q_n:u(t,x)<\sqrt{K}\Lambda^{\alpha(n-1)}M[u]\right\}.$$
	Note that
	$$u_n(t,x)\ge\sup\limits_{(t_n-\delta W(B'_n),t_n]\times 3B'_n}u-u(t,x)\phi_n(x)\ge\sup\limits_{(t_n-\delta W(B'_n),t_n]\times 3B'_n}u-\sqrt{K}\Lambda^{\alpha(n-1)}M[u]=:M_n$$
	on $D_n$, while
	$$M_n\ge u(t_n,y_n)-\sqrt{K}\Lambda^{\alpha(n-1)}M[u]\ge K\Lambda^{\alpha(n-1)}M[u]-\sqrt{K}\Lambda^{\alpha(n-1)}M[u]>0.$$
	According to (\ref{M1}), $(\mathrm{PGL}_1)$ holds under our assumptions. Hence
	$$\frac{u_n(t_n,y_n)}{M_n}\ge C_{\mathrm{PGL}}^{-1}\frac{\hat{\mu}(D_n)}{\hat{\mu}(Q_n)}=:\pi_n.$$
	In particular, at the point $(t_n,y_n)$ we have
	\begin{equation}\label{est-n}
		K\Lambda^{\alpha(n-1)}M[u]\le u(t_n,y_n)\le\sqrt{K}\Lambda^{\alpha(n-1)}\pi_nM[u]+(1-\pi_n)\sup\limits_{(t_n-\delta W(B'_n),t_n]\times 3B'_n}u+A_n(t_n).
	\end{equation}
	
	(4) We estimate $A_n(t_n)$ with proper $\delta''$. Actually, by taking $\delta''$ as in Proposition \ref{tail} we have
	\begin{equation}\label{A_n}
		A_n(t_n)\le D'_1(\delta,\delta'')\mathop{\mathrm{einf}}\limits_{(t_n+2\delta W(B_n),t_n+3\delta W(B_n)]\times\frac{1}{2}B_n}u\le D'_1(\delta,\delta'')\inf\limits_{\frac{1}{2}B_n}u(t''_n,\cdot),
	\end{equation}
	where $t''_n:=t_n+3\delta W(B_n)$. Further, recall that $r\le(6C_W)^{-1/\beta_1}R_0$ by (\ref{r-1}). Hence if we take
	$$\lambda_0:=\frac{1}{12}\wedge\left(\frac{\delta}{2C_W}\right)^{\frac{1}{\beta_1}}\left(\frac{\delta_L}{3}\right)^{\frac{\beta_2}{\beta_1}},$$
	then
	\begin{align*}
		t-t_n<&\ 4\delta W(B)-0=4\delta W(x_0,R_0)\le W\left(y_n,\delta_L\cdot3^{-1}R_0\right);\\
		d(x,y_n)\le&\ d(x,x_0)+d(x_0,y_n)<3\lambda_0R_0\le\delta_LW^{-1}\left(y_n,\frac{\delta}{2}W(x_0,R_0)\right)\\
		\le&\ \delta_LW^{-1}\left(y_n,t-t_n-\frac{\delta}{2C_W}\left(\frac{R_0}{r}\right)^{\beta_1}W(y_0,r)\right)\le\delta_LW^{-1}\left(y_n,t-t_n-3\delta W(y_0,r'_n)\right)=\delta_LW^{-1}(y_n,t''_n).
	\end{align*}
	Therefore, by Corollary \ref{aux-comp} and $(\mathrm{LLE})$,
	\begin{align}
		M[u]=&\ \mathop{\mathrm{einf}}\limits_{Q_+}u\ge\inf\limits_{\frac{1}{2}B_n}u\left(t''_n,\cdot\right)\cdot P_{t-t''_n}^{B(y_n,\frac{1}{3}R_0)}1_{\frac{1}{2}B_n}(x)=\inf\limits_{\frac{1}{2}B_n}u\left(t''_n,\cdot\right)\int_{\frac{1}{2}B_n}p_{t-t''_n}^{B(y_n,\frac{1}{3}R_0)}(x,z)d\mu(z)\notag\\
		\ge&\ \frac{c_L\mu\left(\frac{1}{2}B_n\right)}{V\big(y_n,W^{-1}(y_n,t-t''_n)\big)}\inf\limits_{\frac{1}{2}B_n}u\left(t''_n,\cdot\right)\ge\frac{c_LV\left(y_n,\delta'r'_n\right)}{V\big(y_n,W^{-1}(y_n,4\delta W(B))\big)}\inf\limits_{\frac{1}{2}B_n}u\left(t''_n,\cdot\right)\notag\\
		\ge&\ \frac{1}{D_2(\delta,\delta')}\frac{V(y_n,r'_n)}{V(y_n,R_0)}\inf\limits_{\frac{1}{2}B_n}u\left(t''_n,\cdot\right)\ge\frac{\left(\Lambda^{1-n}\lambda_H\right)^\alpha}{D'_2(\delta,\delta')}\inf\limits_{\frac{1}{2}B_n}u\left(t''_n,\cdot\right).\notag
	\end{align}
	Inserting this inequality into (\ref{A_n}), we obtain
	\begin{equation}\label{tail-n}
		A_n(t)\le\lambda_H^{-\alpha}D'_1(\delta,\delta'')D'_2(\delta,\delta')\Lambda^{\alpha(n-1)}M[u].
	\end{equation}
	
	(5) Now we estimate $\pi_n$. Applying Corollary \ref{aux-comp}, $(\mathrm{PGL}_1)$ and $(\mathrm{LLE})$ successively, we obtain
	\begin{align*}
		\frac{1}{\sqrt{K}\Lambda^{\alpha(n-1)}}\ge&\ \frac{1}{\sqrt{K}\Lambda^{\alpha(n-1)}M[u]}\mathop{\mathrm{einf}}\limits_{Q_+}u\ge\frac{1}{\sqrt{K}\Lambda^{\alpha(n-1)}M[u]}\inf\limits_{\frac{1}{2}\delta''B_n}u(t_n,\cdot)\inf\limits_{(t,x)\in Q_+}P_{t-t_n}^{B\left(y_n,\frac{1}{3}R_0\right)}1_{\frac{1}{2}\delta''B_n}(x)\\
		\ge&\ \frac{C^{-1}_{\mathrm{PGL}}\hat{\mu}(Q_n\setminus D_n)}{\hat{\mu}(Q_n)}\cdot\frac{c_LV(y_n,\frac{1}{2}\delta''r_n)}{V\left(y_n,W^{-1}(y_n,2\delta W(B))\right)}\ge\frac{\hat{\mu}(Q_n\setminus D_n)}{D_3(\delta,\delta'')\frac{\delta}{2}W(\delta''B_n)\mu(B)},
	\end{align*}
	where we applied $(\mathrm{LLE})$ to $p_{t-t_n}^{B\left(y_n,\frac{1}{3}R_0\right)}(x,\cdot)$ in the same way as in Step (4). Consequently,
	$$\frac{\hat{\mu}(D_n)}{\hat{\mu}(Q_n)}=1-\frac{\hat{\mu}(Q_n\setminus D_n)}{\hat{\mu}(Q_n)}\ge 1-\frac{D_3(\delta,\delta'')\frac{\delta}{2}W(\delta''B_n)\mu(B)}{\sqrt{K}\Lambda^{\alpha(n-1)}\frac{\delta}{2}W(\delta''B_n)V(y_n,\delta''B_n)}\ge 1-\frac{D'_3(\delta,\delta',\delta'')}{\lambda_H^\alpha\sqrt{K}}.$$
	Therefore,
	\begin{equation}\label{K-}
		\pi_n\ge\left(2C_{\mathrm{PGL}}\right)^{-1}\quad\quad\mbox{provided that}\quad\quad\sqrt{K}\ge 2D'_3(\delta,\delta',\delta'')\lambda_H^{-\alpha}.
	\end{equation}
	
	(6) Inserting (\ref{tail-n}) and (\ref{K-}) into (\ref{est-n}), we obtain
	$$\sup\limits_{(t_n-\delta W(B'_n),t_n]\times 3B'_n}u\ge\frac{K\Lambda^{\alpha(n-1)}M[u]-\left(\frac{1}{D_4}+\frac{D'_1D'_2}{2D'_3}\right)\sqrt{K}\Lambda^{\alpha(n-1)}M[u]}{1-\frac{1}{D_4}}=K\Lambda^{\alpha(n-1)}M[u]\cdot\frac{1-\frac{1}{\sqrt{K}}\left(\frac{1}{D_4}+\frac{D'_1D'_2}{2D'_3}\right)}{1-\frac{1}{D_4}}.$$
	Therefore, with
	\begin{equation}\label{K+}
		\Lambda=\left(\frac{2D_4-1}{2(D_4-1)}\right)^{1/\alpha}\quad\quad\mbox{and}\quad\quad\sqrt{K}\ge\frac{\frac{1}{D_4}+\frac{D'_1D'_2}{2D'_3}}{1-\Lambda^\alpha\left(1-\frac{1}{D_4}\right)}=\frac{2D_4}{2D_4-1}\left(\frac{1}{D_4}+\frac{D'_1D'_2}{2D'_3}\right),
	\end{equation}
	it follows that
	$$\sup\limits_{(t_n-\delta W(B'_n),t_n]\times 3B'_n}u\ge K\Lambda^{\alpha n}M[u].$$
	Thus by fixing $K>2$ that satisfies (\ref{K-}) and (\ref{K+}), we see there exists $(t_{n+1},y_{n+1})\in(t_n-\delta W(B'_n),t_n]\times 3B'_n$ such that $u(t_{n+1},y_{n+1})\ge K\Lambda^{\alpha n}M[u]$. Hence the proof is completed.
\end{proof}

\subsection{The Moser approach}

\begin{lemma}\label{ep}
	Under conditions $(\mathrm{VD})$, $(\mathrm{FK})$, $(\mathrm{Gcap})$ and $(\mathrm{TJ})$, there exists a constant $C>0$ such that for every three concentric balls $B_0=B(x_0,R)$, $B=B(x_0,R+r)$ and $\Omega=B(x_0,R')$ with $0<R<R+r<R'<\overline{R}$, there is $\phi\in\mathrm{cutoff}(B_0,B)$ such that for all $u\in\mathcal{F}'(\Omega)$ that is bounded in a strict neighbourhood of $\Omega$,
	\begin{equation}\label{ep'}
		\mathcal{E}(u\phi)\le\frac{3}{2}\mathcal{E}(u,u\phi^2)+\frac{C}{\inf\limits_{x\in\Omega}W(x,r)}\int_\Omega u^2d\mu+3\int_{\Omega\times\Omega^c}u(x)u(y)\phi^2(x)dj(x,y).
	\end{equation}
\end{lemma}

\begin{proof}
	By \cite[Remark 5.2]{ghh+2}, the desired result holds when $u\in\mathcal{F}\cap L^\infty$. Now we turn to a general $u\in\mathcal{F}'(\Omega)$.
	
	By definition, there exist $\Omega''\Supset\Omega$ and $w'\in\mathcal{F}$ such that $u=w'$ on $\Omega''$. By the assumption here, $u$ is bouned on some $\Omega'\Supset\Omega$. Take an open set $\Omega'''$ such that $\Omega\Subset\Omega'''\Subset\Omega''\cap\Omega'$, and let $w=w'\phi'$ with some $\phi'\in\mathrm{cutoff}(\Omega''',\Omega')$. Then, $w\in\mathcal{F}\cap L^\infty$ and $u=w$ on $\Omega'''$.
	
	Since $\phi\in\mathrm{cutoff}(B_0,B)\subset\mathcal{F}(B)$, then $u\phi^k=w\phi^k\in\mathcal{F}(B)$ with $k=1,2$, while $\mathrm{supp}(w-u)$ is strictly away from $\overline{\Omega}$. Therefore,
	$$\mathcal{E}\left(w-u,u\phi^2\right)=\mathcal{E}^{(J)}\left(w-u,u\phi^2\right)=-2\int_{\Omega\times\Omega^c}u(x)\phi^2(x)(w-u)(y)dj(x,y).$$
	Combining (\ref{ep'}) applied to $w$, we obtain
	\begin{align*}
		\mathcal{E}(u\phi)=&\ \mathcal{E}(w\phi)\le\frac{3}{2}\mathcal{E}\left(w,w\phi^2\right)+\frac{C}{\inf\limits_{x\in\Omega}W(x,r)}\int_\Omega w^2d\mu+3\int_{\Omega\times\Omega^c}w(x)w(y)\phi^2(x)dj(x,y)\\
		=&\ \frac{3}{2}\mathcal{E}\left(u,u\phi^2\right)+\frac{3}{2}\mathcal{E}\left(w-u,u\phi^2\right)+\frac{C}{\inf\limits_{x\in\Omega}W(x,r)}\int_\Omega u^2d\mu+3\int_{\Omega\times\Omega^c}u(x)\big(u(y)+(w-u)(y)\big)\phi^2(x)dj\\
		=&\ \frac{3}{2}\mathcal{E}(u,u\phi^2)+\frac{C}{\inf\limits_{x\in\Omega}W(x,r)}\int_\Omega u^2d\mu+3\int_{\Omega\times\Omega^c}u(x)u(y)\phi^2(x)dj(x,y),
	\end{align*}
	that is, (\ref{ep'}) holds for $u$.
\end{proof}

\begin{lemma}\label{pmv2}
	Under condition $(\mathrm{VD})$, $(\mathrm{FK})+(\mathrm{Gcap})+(\mathrm{TJ})\Rightarrow(\mathrm{PMV}^-)$.
\end{lemma}

\begin{proof}
	Repeat the proof of \cite[Lemma 5.5]{ghh+2}, but replacing \cite[Implication (5.3)]{ghh+2} (whenever it appears) by Lemma \ref{ep}.
\end{proof}

\begin{proposition}\label{pmv1}
	$(\mathrm{LLE})+(\mathrm{PMV}^-)+(\mathrm{UJS})\Rightarrow(\mathrm{PMV}^+)$.
\end{proposition}

\begin{proof}
	Replacing $W$ by $\delta W$ if necessary, we assume $\delta =1$ in Proposition \ref{tail}. Since $(\mathrm{PMV}^-)$ holds, then the inequality (\ref{pmv2p}) is valid for every target $u$ and ball $B=B(x_0,R)$. By optimizing the right side over $\varepsilon$, we see
	\begin{equation}\label{pmv2pt}
		\mathop{\mathrm{esup}}\limits_{(t_0-\frac{1}{4}W(B),t_0]\times\frac{1}{2}B}u\le C\left(\fint_{Q'}u^2d\hat{\mu}\right)^{\theta/2}\left(\left(\fint_{Q'}u^2d\hat{\mu}\right)^{1/2}\vee\sup\limits_{t_0-\frac{1}{2}W(B)<s\le t_0}\|u_+(s,\cdot)\|_{L^\infty((\frac{1}{2}B)^c)}\right)^{1-\theta},
	\end{equation}
	where $Q':=(t_0-\frac{1}{2}W(B),t_0]\times B$.
	
	With some $0<c_0<\frac{1}{2}$ to be determined later, we denote
	$$B'_0:=c_0B\quad\quad\mbox{and}\quad\quad D_0:=\left(t_0-\frac{1}{2}W(B'_0),t_0+\frac{1}{4}W(B'_0)\right]\times B'_0.$$
	Construct $D_n$ for all $n\ge 1$ inductively as follows: select first $(s_n,x_n)\in D_{n-1}$ such that
	$$\mathop{\mathrm{esup}}\limits_{D_{n-1}}u=\mathop{\mathrm{esup}}\limits_{D_{n-1}\cap Q_n^+}u,$$
	where
	$$Q_n^+:=\left(s_n-\frac{1}{4}W\left(2^{-1}\delta'\delta''B'_n\right),s_n\right]\times\frac{1}{4}\delta'\delta''B'_n\quad\subset\quad\left(s_n-W(B'_n),s_n\right]\times B'_n=:\tilde{Q}_n$$
	with $B'_n=B(x_n,2^{-n}c_0R)$ and $\delta',\delta''$ from Proposition \ref{tail}. Define
	$$D_n=\left(t_0-\frac{1}{2}\sum\limits_{k=0}^nW(B'_k),t_0+\frac{1}{4}\sum\limits_{k=0}^nW(B'_k)\right]\times(2-2^{-k})B'_0.$$
	Similar as Step (2) in the proof of Proposition \ref{main'}, we see there is a sufficiently small constant $c_0$ such that
	$$Q'':=\left(t_0-\frac{1}{2}W(B),t_0+\frac{1}{4}W(B)\right]\times B\quad\supset\quad D_n\supset\tilde{Q}_n\supset Q_n^+\quad\quad\mbox{for all}\quad\quad n\ge 1.$$
	
	Fix an arbitrary $\phi_n\in\mathrm{cutoff}(\delta'B'_n,B'_n)$, and set
	$$Q_n:=\left(s_n-W\left(\delta'\delta''B'_n\right),s_n\right]\times\delta'\delta''B'_n\quad\mbox{and}\quad Q'_n:=\left(s_n+2W(2^{-1}B'_n),s_n+3W(2^{-1}B'_n)\right]\times\frac{1}{2}\delta'B'_n.$$
	By Lemma \ref{trunc} we see
	$$u_n(t,x):=u(t,x)\phi_n(x)+2T_{(t,s_n]\times B'_n}^{(\delta'\delta'')}(u_+)$$
	is globally non-negative and subcaloric on $Q_n$. Further, the second term of $u_n$ is trivially controlled by
	$$T_n:=2T_{(s_n-W(\delta'B'_n),s_n]\times B'_n}^{(\delta'\delta'';1)}(u_+)\le 2C\mathop{\mathrm{einf}}\limits_{Q'_n}u\le 2C\fint_{Q'_n}ud\hat{\mu}.$$
	
	Applying (\ref{pmv2pt}) to $u_n$ on $Q_n$, we obtain
	\begin{equation}\label{mn1}
		M_n:=\mathop{\mathrm{esup}}_{D_{n-1}}u=\mathop{\mathrm{esup}}_{D_{n-1}\cap Q_n^+}u\le\mathop{\mathrm{esup}}_{Q_n^+}u_n\le CA_n^\theta(A_n\vee F_n)^{1-\theta},
	\end{equation}
	where
	$$A_n=\left(\fint_{Q_n}u_n^2d\hat{\mu}\right)^{1/2}\quad\quad\mbox{and}\quad\quad F_n=\sup\limits_{s_n-\frac{1}{2}W(\delta'\delta''B'_n)<s\le s_n}\|u_n(s,\cdot)\|_{L^\infty((\frac{1}{2}\delta'\delta''B'_n)^c)}.$$
	
	Now we estimate $A_n$ and $F_n$ separately. For simplicity, we denote $K:=\fint_{Q''}ud\hat{\mu}$.
	
	On one hand, using the facts $Q_n\subset D_n$ and $d(x_n,x_0)\le 2c_0R$, we see
	\begin{align*}
		A_n^2=&\ \fint_{Q_n}u_n^2d\hat{\mu}\le\fint_{Q_n}u_nd\hat{\mu}\cdot \mathop{\mathrm{esup}}_{D_n}u_n\le\fint_{Q_n}(u+T_n)d\hat{\mu}\cdot(M_{n+1}+T_n)\\
		\le&\ \left(\frac{1}{\hat{\mu}(Q_n)}\int_{Q''}ud\hat{\mu}+\frac{2C}{\hat{\mu}(Q'_n)}\int_{Q''}ud\hat{\mu}\right)\left(M_{n+1}+\frac{2C}{\hat{\mu}(Q'_n)}\int_{Q''}ud\hat{\mu}\right)\le C'2^{2n(\alpha+\beta_2)}(M_{n+1}+K)K;
	\end{align*}
	on the other hand,
	$$F_n=\sup\limits_{s_n-\frac{1}{2}W(\delta'\delta''B'_n)<s\le s_n}\|u_n(s,\cdot)\|_{L^\infty((\frac{1}{2}\delta'\delta''B'_n)^c)}\le\sup\limits_{s_n-\frac{1}{2}W(\delta'\delta''B'_n)<s\le s_n}\|u(s,\cdot)\phi_n\|_{L^\infty}+T_n\le M_{n+1}+C'2^{n(\alpha+\beta_2)}K.$$
	Combining these inequalities, we see
	\begin{align*}
		M_n\le&\ C\left(C'2^{2n(\alpha+\beta_2)}(M_{n+1}+K)K\right)^{\theta/2}\left\{\left(C'2^{2n(\alpha+\beta_2)}(M_{n+1}+K)K\right)^{1/2}\vee\left(M_{n+1}+C'2^{n(\alpha+\beta_2)}K\right)\right\}^{1-\theta}\notag\\
		\le&\ C''2^{n(\alpha+\beta_2)}K^{\theta/2}(M_{n+1}+K)^{1-\theta/2}.
	\end{align*}
	Following the same iteration as in \cite[Proposition 4.2]{ab}, we obtain $M_1\le CK$, which is exactly (\ref{p1mv-p}).
\end{proof}

\begin{theorem}\label{main''}
	Assume that $(\mathrm{wPH}_1)$ holds. Then
	$$(\mathrm{UJS})\Leftrightarrow(\mathrm{PMV}^+)\Leftrightarrow(\mathrm{PHI}^0).$$
\end{theorem}

\begin{proof}
	(1) The implication $(\mathrm{PHI}^0)\Rightarrow(\mathrm{UJS})$ is proved in Proposition \ref{ujs}.
	
	(2) Assume that $(\mathrm{UJS})$ holds. \cite[Theorem 2.2]{L0} we see $(\mathrm{FK})$, $(\mathrm{Gcap})$ (and hence $(\mathrm{cap}_\le)$) hold. Combining Lemma \ref{L2-J}, we yield $(\mathrm{TJ})$. Therefore, by Lemma \ref{pmv2} and then Proposition \ref{pmv1} we obtain $(\mathrm{PMV}^+)$.
	
	(3) Assume that $(\mathrm{PMV}^+)$ holds. Denote the parameters $\delta$ in $(\mathrm{wPH}_1)$ and $(\mathrm{PMV}^+)$ respectively by $\delta_w$ and $\delta_v$. Fix $c_0$ from $(\mathrm{PMV}^+)$ and $\lambda_0$ from $(\mathrm{wPH}_1)$. Fix a ball $B=B(x_0,R)$. Let $\lambda\in(0,c_0\wedge\lambda_0]$ be a constant so small that
	$$r:=W^{-1}\left(x_0,C_Wc_0^{-\beta_2}K_\delta W(\lambda B)\right)\le\lambda_0R,\quad\quad\mbox{where}\quad\quad K_\delta=\frac{3\delta_v}{4\delta_w}\vee 1.$$
	Let
	$$\delta_2=C_Wc_0^{-\beta_2}\left(\frac{\delta_v}{2}+K_\delta\delta_w\right),\quad\quad\delta_1=\delta_2-\frac{\delta_v}{4},\quad\quad\delta_3=3\delta_wC_Wc_0^{-\beta_2}K_\delta\quad\quad\mbox{and}\quad\quad\delta_4=4\delta_wC_Wc_0^{-\beta_2}K_\delta.$$
	We can directly verify that $r\ge c_0^{-1}\lambda R$ and
	$$\delta_wW(x_0,r)\le\delta_2W(\lambda B)-\frac{\delta_V}{2}W(c_0^{-1}\lambda B)<\delta_2W(\lambda B)+\frac{\delta_V}{4}W(c_0^{-1}\lambda B)\le 2\delta_wW(x_0,r).$$
	Hence for every $u$ that is globally non-negative, bounded and caloric on $(0,\delta_4W(\lambda B)]\times B$, we have
	\begin{align*}
		\mathop{\mathrm{esup}}\limits_{(\delta_1W(\lambda B),\delta_2W(\lambda B)]\times\lambda B}u\le&\ \mathop{\mathrm{esup}}\limits_{(\delta_2W(\lambda B)-\frac{\delta_v}{4}W(c_0^{-1}\lambda B),\delta_2W(\lambda B)]\times\lambda B}u\\
		\le&\ C_{\mathrm{PMV}^+}\fint_{(\delta_2W(\lambda B)-\frac{1}{2}\delta_vW(c_0^{-1}\lambda B),\delta_2W(\lambda B)+\frac{1}{4}\delta_vW(c_0^{-1}\lambda B)]\times c_0^{-1}\lambda B}ud\hat{\mu}\\
		\le&\ C'\fint_{(\delta_wW(x_0,r),2\delta_wW(x_0,r)]\times B(x_0,r)}ud\hat{\mu}\\
		\le&\ C'C_{\mathrm{wPH}_1}\mathop{\mathrm{einf}}\limits_{(3\delta_wW(x_0,r),4\delta_wW(x_0,r)]\times B(x_0,r)}u\le C'C_{\mathrm{wPH}_1}\mathop{\mathrm{einf}}\limits_{(\delta_3W(\lambda B),\delta_4W(\lambda B)]\times\lambda B}u.
	\end{align*}
	That is, $(\mathrm{PHI}^0)$ holds with $\delta_i$ given as above, $\lambda_0=\lambda_1=\lambda$ and $C=C'C_{\mathrm{wPH}_1}$.
\end{proof}

There is a further note about $(\mathrm{PMV}^+)$:

\begin{lemma}
	Under condition $(\mathrm{PMV}^+)$, the global heat kernel $p_t(x,y)$ exists, and satisfies the on-diagonal upper estimate (for short, DUE): there exists $C>0$ such that for $\mu$-a.e.\ $x,y\in M$ and all $0<t<\delta W(x,\overline{R})$,
	$$p_t(x,y)\le\frac{C}{V\left(x,W^{-1}(x,t)\right)}.$$
	In particular, $(\mathrm{PMV}^+)\Rightarrow(\mathrm{Nash}^-)$.
\end{lemma}

Combining (\ref{M1}) and Theorem \ref{main''}, it follows directly that
\begin{equation}\label{pmv'nas}
	(\mathrm{PMV}^+)+(\mathrm{wPH})\Rightarrow(\mathrm{PHI}^0).
\end{equation}

\begin{proof}
	Fix an arbitrary $t\in(0,\delta W(x,\overline{R}))$ and set $R:=W^{-1}(x,\delta^{-1}t)$. Then we have
	\begin{equation}\label{pr'due}
		\mathop{\mathrm{esup}}\limits_{\frac{1}{2}B(x,R)}P_tf\le C\fint_{(\frac{1}{2}t,\frac{5}{4}t]\times B(x,R)}P_tfd\hat{\mu}\le C\fint_{\frac{1}{2}t}^{\frac{5}{4}t}\frac{1}{\sqrt{V(x,R)}}\cdot\|P_sf\|_{L^2(B(x,R))}ds\le\frac{C'\|f\|_{L^2(M)}}{\sqrt{V(x,W^{-1}(x,t))}},
	\end{equation}
	where we use H\"older inequality for the second ``$\le$'' and the contractiveness of $P_t$ for the third.
	
	We obtain $(\mathrm{DUE})$ from (\ref{pr'due}) in the same way as in the proof of \cite[Lemma 6.6]{ghh+1}.
	
	Finally, $(\mathrm{DUE})\Rightarrow(\mathrm{Nash}^-)$ obviously since $p_t^B(x,y)\le p_t(x,y)$, implying the last assertion.
\end{proof}

\subsection*{Acknowledgement}

The author would like to thank Professors Jiaxin Hu, Alexander Grigor'yan, Moritz Kassmann and Marvin Weidner for discussion on this research topic. The author ia also grateful to Dr. Lorenzo Portinale for bringing up with many new observations (for example, whether the existence of a killing part is possible).

\end{document}